\documentclass{article}

\usepackage{arxiv}

\usepackage[utf8]{inputenc} 
\usepackage[T1]{fontenc}    
\usepackage{hyperref}       
\usepackage{url}            
\usepackage{booktabs}       
\usepackage{amsfonts}       
\usepackage{amsmath}
\usepackage{amssymb}
\usepackage{nicefrac}       
\usepackage{microtype}      
\usepackage{graphicx}
\usepackage{hyperref}
\usepackage{multirow}
\usepackage{tikz}
\usepackage{cancel}
\usepackage{amsmath}
\usepackage{amsthm}
\usepackage{pgfplots}
\usepackage{longtable}
\usepackage{float}
\usepackage{algorithm}
\usepackage{algpseudocode}

\usetikzlibrary{positioning}
\usetikzlibrary{decorations.pathreplacing,calligraphy}
\usetikzlibrary{intersections}
\usepgfplotslibrary{fillbetween}

\newtheorem{lemma}{Lemma}
\newtheorem{definition}{Definition}
\newtheorem{observation}{Observation}
\graphicspath{ {./} }

\providecommand{\citet}{\cite}
\providecommand{\citep}{\cite}
\providecommand{\citealp}{\cite}

\title{A parallel pull labelling algorithm for the resource constrained shortest path problem}

\author{
  Bjørn Petersen \\
  Flowty ApS\\
  Denmark\\
  \texttt{bjorn@flowty.ai} \\
\And
 Simon Spoorendonk \\
  Flowty ApS\\
  Denmark\\
  \texttt{simon@flowty.ai} \\
}

\begin{document}
\maketitle
\begin{abstract}
The Resource Constrained Shortest Path Problem (RCSPP) is a fundamental combinatorial optimisation problem in
which the goal is to find a least-cost path in a directed graph subject to one or more resource constraints.
\emph{Pull} labelling, where a vertex gathers labels from its predecessors rather than pushing them to its successors,
is a classical idea; we turn it into a parallel algorithm whose immutable, contention-free bucket storage is designed to scale on modern multi-core hardware.
Its central component is an acyclic \emph{dependency bucket graph} that orders the creation of labels so that buckets
can be processed concurrently, without conflicting accesses, and stored as immutable objects. On top of this we introduce
i) a highly parallelisable approach at the label-bucket level,
ii) an extension to bi-directional search with a dynamic midpoint that emerges from the bucket processing order, and iii) a vectorised dominance criterion
that uses vector instructions to speed-up the label comparison with another level of parallelisation.
Compared to a baseline version of the algorithm the optimisations result in a speed-up of about 18 times
on a set of hard instances and up to 274 times on the instance with the largest speed-up.
Against an open implementation of the state-of-the-art bucket graph labelling algorithm, run on the same hardware and instances with each solver in its best parallel configuration,
the pull algorithm is 1.9 to 2.4 times faster.
The proposed algorithm demonstrates significant computational improvements that may enhance
the efficiency of column generation frameworks incorporating resource constrained shortest path
sub-problems, potentially enabling the efficient solution of larger-scale instances in routing, scheduling, supply chain
and transportation network optimisation applications.

\end{abstract}

\keywords{Resource constrained shortest path problem, Labelling algorithm, Parallelisation}

\newif\ifijoc\ijocfalse

\section{Introduction}
\label{sec:introduction}

The Resource Constrained Shortest Path Problem (RCSPP) is a fundamental combinatorial optimisation problem in 
which the goal is to find a least-cost path in a directed graph subject to one or more resource constraints \citep{irnich2005, irnich2008}.
Formally, given a directed graph $G=(V,E)$ with vertices $V$, edges $E$, a source $s \in V$, a sink $t \in V$, 
and a set of resources $R$. Each edge $e \in E$ has an associated cost $c_e$ and resource consumptions $q^r_e$ for each $r \in R$. 
Additionally, for each resource $r \in R$, bounds $[a^r_v, b^r_v]$ are associated with each vertex $v \in V$ or $[a^r_e, b^r_e]$ for each edge $e \in E$.

Each resource $r$ has a resource extension function (REF) $f_r$ giving its evolution across an edge $(i,j)$, $s^r_j = f_r(s^r_i, q^r_{ij})$; examples
are time-window propagation $s^r_j = \max\{a^r_j, s^r_i + q^r_{ij}\}$ and capacity accumulation $s^r_j = s^r_i + q^r_{ij}$. A REF may also depend on the
edge ($f_{r,e}$), allowing resets such as daily driving time in multi-day routing; such non-monotone resets are permitted by the problem but fall
outside the sufficient dominance condition of \autoref{sec:vectorised} (\autoref{lemmaSufficientCondition}) and are not exercised in our experiments.
A state $S_i = (s^r_i)_{r\in R}$ is feasible if all resources satisfy their bounds, e.g., $a^r_v \le s^r_v \le b^r_v$.
Binary or bit-vector resources can be used to encode visit patterns like pickup and delivery or enforce elementary paths
\citep{feillet2004}. 

The RCSPP is NP-hard and arises as a core sub-problem in many large-scale applications, both as a standalone 
problem and as a pricing sub-problem within column generation and branch-and-price algorithms 
\citep{desaulniers2005, desrosiers2024, uchoa2024}.

The RCSPP is commonly solved via labelling algorithms, a dynamic programming method
that extend partial paths one edge at a time while maintaining feasibility and dominance conditions \citep{irnich2005}. 
A label represents a partial path and stores the accumulated cost and resource state,
while dominated labels, those that are no better in any dimension (cost and every resource) than some other label at the same vertex, are pruned to limit the search space.

Push labelling algorithms (e.g., \citealp{irnich2005, baldacci2011, sadykov2020}) handle dominance and feasibility pruning
but can be memory intensive due to frequent dynamic insertions and removals
of labels into the label storage of the vertices extended to. The algorithm parallelises well to two threads (forward and backward search) 
but more threads can lead to higher contention for the label storage and slower performance. 

Bi-directional search was introduced by \citet{salani2006} where a deterministically chosen midpoint was used. 
Later \citep{tilk2017} introduced a dynamic midpoint selection strategy that balances the number of forward and backward labels, and enhancements
to bi-directional dynamic programming for the RCSPP have continued to be studied since \citep{salani2024}.
We follow this line of research and extend the pull labelling algorithm to a bi-directional algorithm with dynamic midpoint selection.

\citep{sadykov2020} uses buckets to group labels into intervals of resource accumulation within a vertex, 
this allows for the optimisation of dominance checks by only comparing with labels in relevant buckets 
and use quick rejections on minimal label costs for a bucket.

When elementary paths are required, the $ng$-path relaxation of \citet{baldacci2011} is common: each vertex $i$ has a neighbourhood $ng_i \subseteq V$,
and a vertex $j$ may be revisited only once the path has left the neighbourhoods of all vertices visited in between (i.e.\ $j$ may not reappear while
$j \in ng_k$ for some still-``remembered'' $k$). Smaller neighbourhoods give a weaker but cheaper relaxation; they may vary by vertex and be grown
dynamically up to a maximum size \citep{bulhoes2018}, the visit-memory carried as a bit-vector resource.

In a column generation context where master problems are set partitioning problems cutting planes directly 
on master problem variables are effective to find good relaxation \citep{desrosiers2024, uchoa2024}. 
For the RCSPP this gives rise to a non-linear cost function that was handled elegantly by \citet{jepsen2008} and
later improved by \citet{pecin2017} using a limited memory technique. Such techniques are combined in the state-of-the-art exact
branch-cut-and-price algorithms for the VRPTW \citep{pecin2017b}.

The pull principle itself is not new. Already \citet{desrosiers1995}, in the context of time constrained routing and scheduling,
organise labels into buckets and compute a vertex's labels from those of its predecessors -- the essential pull idea. Parallel and vectorised labelling
has also been studied outside the column generation setting: \citet{lu2021} accelerate \emph{exact} constrained shortest paths on GPUs, organising
labels into resource buckets with dominance pruning and Single Instruction Multiple Data (SIMD) execution. That line of work, like the bucket graph of \citet{sadykov2020}, still
maintains a mutable label store from which dominated labels are removed. What we add is a \emph{pull} formulation, in the column generation pricing
setting with an $ng$-relaxation, in which an acyclic dependency bucket graph makes the label store \emph{immutable}: labels, once created, are never
revisited or removed. Because each bucket is pulled to independently of the others, its dominance checks run in parallel without conflicting memory
accesses, and because the acyclic graph (\autoref{sec:mono}) fixes the creation order, labels are created already dominance-ordered; this is what makes
the parallelisation, the bi-directional extension, and the vectorised dominance check work together. This is the gap we address.

The main contribution of this paper is this parallel pull labelling algorithm whose acyclic dependency bucket graph makes label storage immutable, so
that label-buckets are processed concurrently without coordinating access while the dominance machinery is preserved. On this basis we further
contribute a bi-directional search whose dynamically determined midpoint emerges directly from the bucket processing order rather than being fixed a
priori, and a vectorised dominance criterion that exploits vector instructions (SIMD) to accelerate label comparisons. None of the three components is
novel in isolation: pull labelling, bi-directional search, and SIMD are all known, and prior work has parallelised and vectorised labelling
\citep{lu2021}. What is new, to our knowledge, is combining them in a single immutable-storage pull scheme with contention-free bucket access for
column generation pricing.

The remainder of this paper is organised as follows.
\autoref{sec:preliminaries} reviews preliminaries on push and pull labelling algorithms, 
\autoref{sec:mono} presents the parallel mono-directional pull labelling algorithm,
\autoref{sec:bidirectional} extends the algorithm to a bi-directional version with dynamic midpoint selection, 
\autoref{sec:vectorised} introduces the vectorised dominance criterion, and \autoref{sec:experiments} reports 
 computational experiments demonstrating the effectiveness of the proposed algorithm. Finally, \autoref{sec:conclusion} concludes the paper.

\section{Pull and Push Labelling Algorithms}
\label{sec:preliminaries}

The basic component in labelling algorithms is the label. A label is an object representing a (partial) path and contains properties
like the pointer to the parent label, the vertex, the cost, and the state with current resource accumulations.

The first label at the origin represents an empty path with initial values set accordingly, usually at zero cost and minimum resource accumulations. 
For bi-directional labelling an additional first label is created at the target vertex with initial values set 
at minimum cost and maximum resource accumulations. 
REFs are reversed for backward labelling such that resource accumulation in general is decreasing when going backwards. 

Labelling is an enumeration scheme where paths are grown by extending labels one edge at a time. 
The extension function propagates the path over an edge by applying the REFs to the resource accumulations. To guarantee that the
labelling terminates, the buckets are ordered by a designated \emph{monotone} resource: a resource whose accumulation is non-decreasing
along every edge of the graph (e.g.\ elapsed time or accumulated capacity), so that each extension strictly advances the path in that resource.
This monotone resource need not be the cost, and other resources may behave non-monotonically. For graphs with negative-cost cycles such a
monotone resource is what guarantees that the bucket dependency graph is acyclic and that the search progresses.

\begin{figure}[htbp]
\begin{center}
\begin{minipage}{.4\textwidth}
\begin{center}
\begin{tikzpicture}[
    nodeTiny/.style={rectangle, draw=, fill=white, very thick, minimum width=8mm, minimum height=8mm},
    ]
    \node[nodeTiny](a) at (0.0, 1.2){\Large $a$};
    \node[nodeTiny](b) at (0.0, 0.0){\Large $b$};
    \node[nodeTiny](c) at (0.0, -1.2){\Large $c$};

    \node[nodeTiny](d) at (1.6, 1.2){\Large $d$};
    \node[nodeTiny](e) at (1.6, 0.0){\Large $e$};
    \node[nodeTiny](f) at (1.6, -1.2){\Large $f$};

    \draw[->, line width=0.25mm] (a) -- (e);
    \draw[->, line width=0.25mm] (b) -- (e);
    \draw[->, line width=0.25mm] (c) -- (e);
\end{tikzpicture}
\end{center}
\end{minipage}
\begin{minipage}{.4\textwidth}
\begin{center}
\begin{tikzpicture}[
    nodeTiny/.style={rectangle, draw=, fill=white, very thick, minimum width=8mm, minimum height=8mm},
    ]
    \node[nodeTiny](a) at (0.0, 1.2){\Large $a$};
    \node[nodeTiny](b) at (0.0, 0.0){\Large $b$};
    \node[nodeTiny](c) at (0.0, -1.2){\Large $c$};

    \node[nodeTiny](d) at (1.6, 1.2){\Large $d$};
    \node[nodeTiny](e) at (1.6, 0.0){\Large $e$};
    \node[nodeTiny](f) at (1.6, -1.2){\Large $f$};

    \draw[->, line width=0.25mm] (b) -- (d);
    \draw[->, line width=0.25mm] (b) -- (e);
    \draw[->, line width=0.25mm] (b) -- (f);
\end{tikzpicture}
\end{center}
\end{minipage}
\end{center}
\caption{Pull (left) where labels are pulled into bucket $e$ from buckets $a$, $b$, and $c$. Push (right) where labels are pushed from bucket $e$ to buckets $a$, $b$, and $c$.}
\label{fig:pull-push}
\end{figure}
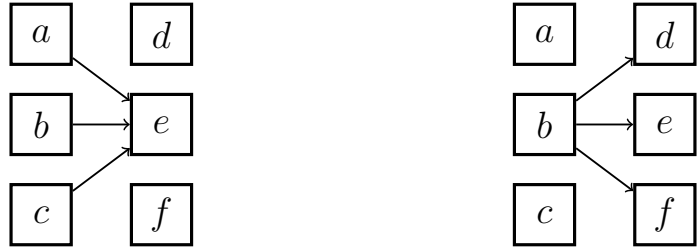

A pull algorithm gathers labels into a vertex from its predecessors, whereas a push algorithm pushes them from a vertex to its successors
(\autoref{fig:pull-push}); in both, extended labels are stored at the vertex they arrive at. Pulling updates a vertex's labels independently of the
other vertices, rather than writing into every neighbour, which is the key advantage for parallelisation: the updates run in parallel without
conflicting writes. It also controls the order in which labels are created, letting us test a new label for dominance before storing or discarding it.
As shown in \autoref{sec:vectorised} (\autoref{lemmaDominanceOrder}), pulling in this controlled order ensures a kept label is never dominated by a
later one, so it never has to be removed and labels can be stored as immutable data structures.

\subsection{Buckets}

Following \citet{sadykov2020}, buckets group the labels of a vertex into intervals of the monotone resource. In a push setting they speed up
dominance by restricting comparisons to relevant buckets and by quick-rejecting on a bucket's minimal cost; in the pull setting these quick rejections
are unnecessary, because the pulling scheme delivers exactly the relevant labels to each bucket. The buckets of a vertex are a scheduling and locality
device, not a partition that hides labels from one another: a new label is checked against all of the vertex's stored labels of no greater monotone
resource, which span the lower buckets too, so cross-bucket dominance within a vertex is fully accounted for. In a pull algorithm the buckets are also
the unit of parallelism.

Buckets induce a graph between buckets; when bucket sizes are small enough this bucket graph is acyclic, which is essential. A large bucket size can
create cycles, so that buckets are processed and pulled several times and labels created out of dominance order later have to be re-dominated,
defeating immutable storage. \citet{sadykov2020} permit small cycles in their push-based bucket algorithm, where the memory impact is smaller because
they re-dominate created labels anyway.

Formally, let $\rho$ denote the monotone resource on which buckets are defined and let $[a^\rho_v, b^\rho_v]$ be its window at vertex $v$.
We set the bucket (step) size of $v$ to the smallest consumption of $\rho$ over the edges incident to $v$,
$$\Delta_v = \min\{\, q^\rho_e : e \in \delta(v) \,\},$$
where $\delta(v)$ are the edges into and out of $v$ that lie on some feasible path. We assume throughout that $\rho$ is strictly
consumed on every edge, $q^\rho_e > 0$, so $\Delta_v > 0$ and $N_v$ is well defined. The number of buckets at $v$ is then
$$N_v = \left\lfloor \frac{b^\rho_v - a^\rho_v}{\Delta_v} \right\rfloor + 1,$$
and a resource value $s \in [a^\rho_v, b^\rho_v]$ falls in bucket $\lfloor (s - a^\rho_v)/\Delta_v \rfloor$.
Choosing $\Delta_v$ as the smallest incident consumption ensures every edge into $v$ advances by at least one of $v$'s buckets; equivalently
$\Delta_v \le q^\rho_e$ for every incident edge $e$. Under this choice the forward dependency graph is acyclic. Order the buckets lexicographically by
$(\rho\text{-interval lower bound},\,\text{bucket index})$. Along every inter-vertex dependency edge $x_m \to y_n$ the destination's lower bound
strictly exceeds the source's: the label's value of $\rho$ advances by $q^\rho_e > 0$ and, because the step $\Delta_y \le q^\rho_e$ does not overshoot a
whole destination bucket, $a^\rho_y + n\,\Delta_y > a^\rho_x + m\,\Delta_x$. Every intra-vertex dependency edge $f^v_{n-1}\to f^v_n$ goes to a strictly
higher bucket index at the same vertex. Hence the key strictly increases on every edge, so the graph cannot contain a cycle. This is precisely why the
step size must not exceed the smallest incident consumption: coarser buckets can violate $\Delta_v \le q^\rho_e$ and reintroduce cycles. Two buckets
may be joined by several edges of the original graph; the dependency counts simply accumulate these, so multiplicity is handled without affecting
acyclicity (see \autoref{sec:mono}).

The bucket size trades off overhead against granularity: larger buckets reduce parallelism overhead, smaller ones give fewer dependencies but sparser
buckets, and the size must stay small enough to avoid cycles. The choice $\Delta_v$ above is the largest size that still guarantees acyclicity.

\section{Parallel Mono-directional Pull Labelling Algorithm}
\label{sec:mono}

At the core of the algorithm is the acyclic dependency graph, which tracks when a bucket may be processed -- that is, when all its dependencies are
done, so that buckets with no unprocessed dependencies can be processed concurrently (processing a bucket is a \emph{job}). The first label is placed in
the origin's first bucket, which has no dependencies; all buckets depending only on it are then ready, initiating the algorithm.

\ifijoc Figure~EC.2 in the online supplement shows\else\autoref{fig:dependency-graph-initial} shows\fi\ the initial state of the dependency graph
between buckets in vertices $x$ and $y$ with size $4$ and $6$, respectively.
An edge from $x$ to $y$ has a resource consumption of $6$ units and from $y$ to $x$ of $4$ units. The edges in the graph represent 
the dependency between the buckets, e.g., $x_2$ depends on $x_1$ and $y_1$, and $y_2$ depends on $x_1$ and $x_2$.

\ifijoc\else
\begin{figure}[htbp]
\begin{center}
\begin{tikzpicture}[
    nodeTiny/.style={rectangle, draw=, fill=white, very thick, minimum width=10mm, minimum height=7.5mm},
    nodeSmall/.style={rectangle, draw=, fill=white, very thick, minimum width=10mm, minimum height=10mm},
    nodeLarge/.style={rectangle, draw=, fill=white, very thick, minimum width=10mm, minimum height=15mm},
    ]

    \node[nodeTiny](x3) at (0.0, -0.125){\LARGE $x_3$};
    \draw[-] (-0.5, 0.0) -- (0.5, 0.0);
    \draw[-] (-0.5, -0.25) -- (0.5, -0.25);

    \node[nodeSmall](x2) at (0.0,-1.0){\huge $x_2$};
    \draw[-] (-0.5, -0.75) -- (0.5, -0.75);
    \draw[-] (-0.5, -1.0) -- (0.5, -1.0);
    \draw[-] (-0.5, -1.25) -- (0.5, -1.25);

    \node[nodeSmall](x1) at (0.0,-2.0){\huge $x_1$};
    \draw[-] (-0.5, -1.75) -- (0.5, -1.75);
    \draw[-] (-0.5, -2.0) -- (0.5, -2.0);
    \draw[-] (-0.5, -2.25) -- (0.5, -2.25);
    \node[draw=none] at (-1.45,-1.0) {size $=4$};
    \draw[decorate, decoration={brace, raise=5pt}, line width=0.25mm] (-0.40,-1.50) -- (-0.40,-0.50);

    \node[draw=none] at (-0.8, 0.25) {$b_x$};
    \node[draw=none] at (-0.8,-2.50) {$a_x$};

    \node[nodeSmall](y4) at (3.0, 2.25){\huge $y_4$};
    \draw[-] (2.5, 2.50) -- (3.5, 2.50);
    \draw[-] (2.5, 2.25) -- (3.5, 2.25);
    \draw[-] (2.5, 2.00) -- (3.5, 2.00);

    \node[nodeLarge](y3) at (3.0, 1.00){\huge $y_3$};
    \draw[-] (2.5, 1.50) -- (3.5, 1.50);
    \draw[-] (2.5, 1.25) -- (3.5, 1.25);
    \draw[-] (2.5, 1.0) -- (3.5, 1.0);
    \draw[-] (2.5, 0.75) -- (3.5, 0.75);
    \draw[-] (2.5, 0.50) -- (3.5, 0.50);

    \node[nodeLarge](y2) at (3.0,-0.50){\huge $y_2$};
    \draw[-] (2.5, 0.00) -- (3.5, 0.00);
    \draw[-] (2.5, -0.25) -- (3.5, -0.25);
    \draw[-] (2.5, -0.50) -- (3.5, -0.50);
    \draw[-] (2.5, -0.75) -- (3.5, -0.75);
    \draw[-] (2.5, -1.00) -- (3.5, -1.00);

    \node[nodeLarge](y1) at (3.0,-2.00){\huge $y_1$};
    \draw[-] (2.5, -1.50) -- (3.5, -1.50);
    \draw[-] (2.5, -1.75) -- (3.5, -1.75);
    \draw[-] (2.5, -2.0) -- (3.5, -2.0);
    \draw[-] (2.5, -2.25) -- (3.5, -2.25);
    \draw[-] (2.5, -2.50) -- (3.5, -2.50);
    \node[draw=none] at (4.45, -0.50) {size $=6$};
    \draw[decorate, decoration={brace, raise=5pt}, line width=0.25mm] (3.40, 0.25) -- (3.40, -1.25);

    \node[draw=none] at (3.8, 2.75) {$b_y$};
    \node[draw=none] at (3.8,-2.75) {$a_y$};

    \draw (0.48, 0.125) -- (2.48, 1.625) node [midway, above, sloped] (TextNode) {6 units};
    \draw[->, line width=0.25mm] (0.48, 0.125) -- (2.48, 1.625);
    \draw[->, line width=0.25mm] (0.48, -0.625) -- (2.48, 0.875);
    \draw[->, line width=0.25mm] (0.48, -1.375) -- (2.48, 0.125);
    \draw[->, line width=0.25mm] (0.48, -2.375) -- (2.48, -0.875);

    \draw[->, line width=0.25mm] (2.48, -1.125) -- (0.48, -0.125);
    \draw (2.48, -1.125) -- (0.48, -0.125) node [midway, above, sloped] (TextNode) {4 units};
    \draw[->, line width=0.25mm] (2.48, -1.375) -- (0.48, -0.375);
    \draw[->, line width=0.25mm] (2.48, -2.375) -- (0.48, -1.375);

\end{tikzpicture}

$y_1 \rightarrow \{x_2, x_3, y_2\}$,
$x_1 \rightarrow \{y_2, x_2\}$,
$x_2 \rightarrow \{y_2, y_3, x_3\}$,
$y_2 \rightarrow \{x_3, y_3\}$,
$x_3 \rightarrow \{y_3\}$,
$y_3 \rightarrow \{y_4\}$
\end{center}

\caption{Initial state of the dependency graph.}
\label{fig:dependency-graph-initial}
\end{figure}
\fi

For a series of edges in the dependency graph a simple observation can be made that can be used to simplify the graph. Here $z_1 \rightarrow \{z_2,\dots\}$
means that bucket $z_1$ must be processed before the listed buckets.

\begin{observation}[Implied edges]
$z_1 \rightarrow \{z_2, z_3 \} \land z_2 \rightarrow \{z_3 \} \Leftrightarrow z_1 \rightarrow \{z_2 \} \land z_2 \rightarrow \{z_3 \}$
\end{observation}

That is, if $z_1$ precedes $z_2$ and $z_2$ precedes $z_3$, the dependency $z_1 \rightarrow z_3$ is already implied by transitivity and need not be
stored. Hence the edges that are implied by the others can be removed, leaving the transitive reduction of the dependency graph.

$y_1 \rightarrow \{x_2, \cancel{x_3}, \cancel{y_2}\}$,
$x_1 \rightarrow \{\cancel{y_2}, x_2\}$,
$x_2 \rightarrow \{y_2, \cancel{y_3}, \cancel{x_3}\}$, \\
$y_2 \rightarrow \{x_3, \cancel{y_3}\}$,
$x_3 \rightarrow \{y_3\}$,
$y_3 \rightarrow \{y_4\}$

This results in the following dependency graph\ifijoc\ (Figure~EC.3 in the online supplement)\else, as depicted in \autoref{fig:dependency-graph-minimal}\fi.
Here $x_1$ and $y_1$ are the two initial buckets with no dependencies; the remaining buckets are then processed in the order
$x_1, y_1 \rightarrow x_2 \rightarrow y_2 \rightarrow x_3 \rightarrow y_3 \rightarrow y_4$, where $x_1$ and $y_1$ are ready first and may be
processed in either order or in parallel.

\ifijoc\else
\begin{figure}[htbp]
\begin{center}
\begin{tikzpicture}[
    nodeTiny/.style={rectangle, draw=, fill=white, very thick, minimum width=10mm, minimum height=7.5mm},
    nodeSmall/.style={rectangle, draw=, fill=white, very thick, minimum width=10mm, minimum height=10mm},
    nodeLarge/.style={rectangle, draw=, fill=white, very thick, minimum width=10mm, minimum height=15mm},
    ]

    \node[nodeTiny](x3) at (0.0, -0.125){\LARGE $x_3$};
    \draw[-] (-0.5, 0.0) -- (0.5, 0.0);
    \draw[-] (-0.5, -0.25) -- (0.5, -0.25);

    \node[nodeSmall](x2) at (0.0,-1.0){\huge $x_2$};
    \draw[-] (-0.5, -0.75) -- (0.5, -0.75);
    \draw[-] (-0.5, -1.0) -- (0.5, -1.0);
    \draw[-] (-0.5, -1.25) -- (0.5, -1.25);

    \node[nodeSmall](x1) at (0.0,-2.0){\huge $x_1$};
    \draw[-] (-0.5, -1.75) -- (0.5, -1.75);
    \draw[-] (-0.5, -2.0) -- (0.5, -2.0);
    \draw[-] (-0.5, -2.25) -- (0.5, -2.25);
    \node[draw=none] at (-1.45,-1.0) {size $=4$};
    \draw[decorate, decoration={brace, raise=5pt}, line width=0.25mm] (-0.40,-1.50) -- (-0.40,-0.50);

    \node[draw=none] at (-0.8, 0.25) {$b_x$};
    \node[draw=none] at (-0.8,-2.50) {$a_x$};

    \node[nodeSmall](y4) at (3.0, 2.25){\huge $y_4$};
    \draw[-] (2.5, 2.50) -- (3.5, 2.50);
    \draw[-] (2.5, 2.25) -- (3.5, 2.25);
    \draw[-] (2.5, 2.00) -- (3.5, 2.00);

    \node[nodeLarge](y3) at (3.0, 1.00){\huge $y_3$};
    \draw[-] (2.5, 1.50) -- (3.5, 1.50);
    \draw[-] (2.5, 1.25) -- (3.5, 1.25);
    \draw[-] (2.5, 1.0) -- (3.5, 1.0);
    \draw[-] (2.5, 0.75) -- (3.5, 0.75);
    \draw[-] (2.5, 0.50) -- (3.5, 0.50);

    \node[nodeLarge](y2) at (3.0,-0.50){\huge $y_2$};
    \draw[-] (2.5, 0.00) -- (3.5, 0.00);
    \draw[-] (2.5, -0.25) -- (3.5, -0.25);
    \draw[-] (2.5, -0.50) -- (3.5, -0.50);
    \draw[-] (2.5, -0.75) -- (3.5, -0.75);
    \draw[-] (2.5, -1.00) -- (3.5, -1.00);

    \node[nodeLarge](y1) at (3.0,-2.00){\huge $y_1$};
    \draw[-] (2.5, -1.50) -- (3.5, -1.50);
    \draw[-] (2.5, -1.75) -- (3.5, -1.75);
    \draw[-] (2.5, -2.0) -- (3.5, -2.0);
    \draw[-] (2.5, -2.25) -- (3.5, -2.25);
    \draw[-] (2.5, -2.50) -- (3.5, -2.50);
    \node[draw=none] at (4.45, -0.50) {size $=6$};
    \draw[decorate, decoration={brace, raise=5pt}, line width=0.25mm] (3.40, 0.25) -- (3.40, -1.25);

    \node[draw=none] at (3.8, 2.75) {$b_y$};
    \node[draw=none] at (3.8,-2.75) {$a_y$};

    \draw (0.48, 0.125) -- (2.48, 1.625) node [midway, above, sloped] (TextNode) {6 units};
    \draw[->, line width=0.25mm] (0.48, 0.125) -- (2.48, 1.625);
    \draw[->, line width=0.25mm] (0.48, -1.375) -- (2.48, 0.125);

    \draw[->, line width=0.25mm] (2.48, -1.125) -- (0.48, -0.125);
    \draw (2.48, -1.125) -- (0.48, -0.125) node [midway, above, sloped] (TextNode) {4 units};
    \draw[->, line width=0.25mm] (2.48, -2.375) -- (0.48, -1.375);

\end{tikzpicture}

$y_1 \rightarrow \{x_2\}$,
$x_1 \rightarrow \{x_2\}$,
$x_2 \rightarrow \{y_2\}$,
$y_2 \rightarrow \{x_3\}$,
$x_3 \rightarrow \{y_3\}$,
$y_3 \rightarrow \{y_4\}$
\end{center}
\caption{Minimal dependency graph.}
\label{fig:dependency-graph-minimal}
\end{figure}
\fi

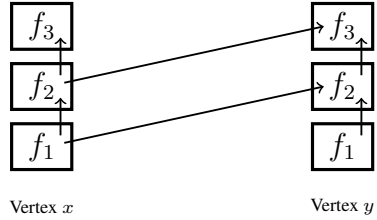
\begin{figure}[htbp]
    \begin{center}
    \colorlet{shadecolor}{gray!34}
    \begin{tikzpicture}[
        nodeWhite/.style={rectangle, draw=, fill=white, very thick, minimum width=8mm, minimum height=6mm},
        nodeGrey/.style={rectangle, draw=, fill=shadecolor, very thick, minimum width=8mm, minimum height=6mm},
        ]
    
        \node[nodeWhite](f3_1) at (0.0, 0.4){\large $f_3$};
        \node[nodeWhite](f2_1) at (0.0, -0.4){\large $f_2$};
        \node[nodeWhite](f1_1) at (0.0, -1.2){\large $f_1$};

        \node[nodeWhite](f3_2) at (4.0, 0.4){\large $f_3$};
        \node[nodeWhite](f2_2) at (4.0, -0.4){\large $f_2$};
        \node[nodeWhite](f1_2) at (4.0, -1.2){\large $f_1$};

        \draw[->, line width=0.25mm] (0.3, -0.35) -- (3.75, 0.4);
        \draw[->, line width=0.25mm] (0.3, -1.15) -- (3.75, -0.4);        
        \draw[->, line width=0.25mm] (0.25, -0.25) -- (0.25, 0.25);
        \draw[->, line width=0.25mm] (0.25, -1.05) -- (0.25, -0.55);
        \draw[->, line width=0.25mm] (4.25, -0.25) -- (4.25, 0.25);
        \draw[->, line width=0.25mm] (4.25, -1.05) -- (4.25, -0.55);

        \node[draw=none] at (0.0, -2.0) {\scriptsize Vertex $x$};
        \node[draw=none] at (4.0, -2.0) {\scriptsize Vertex $y$};
    \end{tikzpicture}
    \end{center}
    \caption{Dependency graph for mono-directional pull labelling algorithm. }
\label{fig:dependency-graph-mono}
\end{figure}

A sketch of the dependency graph for the mono-directional search is shown in \autoref{fig:dependency-graph-mono}. 
The sketch exemplifies that forward buckets depend on forward buckets from other vertices where the labels can be 
pulled from and updated according to the REFs, i.e., the edges illustrated from buckets in vertex $x$ to buckets in vertex $y$. 
Also, forward buckets depend on the forward bucket below, with less resource accumulation, i.e.,
the vertical edges between buckets within the vertex. The illustration depicts equal sized buckets 
in both vertices which is not necessarily the case in practice although the concept remains the same.

A bucket is said to be released when ready to be processed, i.e., when all its dependencies have been processed. We write $f^y_n$ for the $n$-th
forward bucket of vertex $y$ (and later $b^y_n$, $s^y_n$ for backward and splice); the superscript is the vertex and the subscript the bucket index.
This reuse of the letter $f$ is distinct from the REF $f_r$ of \autoref{sec:preliminaries}, which carries a resource subscript $r$; the meaning is
always clear from the subscript.
In general, the forward bucket $f^y_n$ at vertex $y$ depends on the previous forward bucket $f^y_{n-1}$ of the same vertex and,
for every in-edge $(x,y)$, on the forward bucket of $x$ whose labels can still be pulled into bucket $n$ across that edge:

\begin{tabular}{ll}
     $\Call{Release}{f^y_n}$ & \hspace{-0.5cm} $\gets \Call{processed}{f^y_{n-1}} \land \bigwedge_{(x,y)\in E} \Call{processed}{f^x_{\pi(x,y,n)}}$ \\
\end{tabular}

Here $\pi(x,y,n)$ denotes the last forward bucket of $x$ that can feed bucket $n$ of $y$ over edge $(x,y)$, i.e.\ the largest bucket index
$m$ with $a^\rho_x + m\,\Delta_x + q^\rho_{xy} < a^\rho_y + (n+1)\,\Delta_y$; the inequality is strict because bucket $n$ is the half-open
interval $[a^\rho_y + n\,\Delta_y,\, a^\rho_y + (n+1)\,\Delta_y)$. It accounts for the buckets of $x$ and $y$ having different
sizes $\Delta_x, \Delta_y$. If the resulting index is negative, no bucket of $x$ can feed bucket $n$ of $y$ and the cross-edge dependency is omitted;
if it is at least $N_x$, every bucket of $x$ can, and the dependency binds on $x$'s last bucket $N_x-1$. When all buckets share one size and edge
consumptions equal that size this reduces to $\pi(x,y,n)=n-1$, the schematic case of \autoref{fig:dependency-graph-mono}.

The number of dependencies for a bucket is handled implicitly by keeping count of remaining dependencies.
When a bucket is processed all buckets depending on it are decremented and if a dependent bucket reaches a 
dependency count of zero that bucket is released. See \autoref{alg:remove-dependency}.

\begin{algorithm}[H]
\caption{Release buckets when their dependencies have been processed.}
\label{alg:remove-dependency}
\begin{algorithmic}
\Procedure{RemoveDependency}{bucket: $b$}
\ForAll {$i$ blocked by $b$}
    \State $DependencyCount{[i]} \gets DependencyCount{[i]} - 1$
    \If{$DependencyCount{[i]} = 0$}
        \State \Call{ReleaseBucket}{i}
    \EndIf
\EndFor
\EndProcedure
\end{algorithmic}
\end{algorithm}

\subsection{Parallel Read/Write}
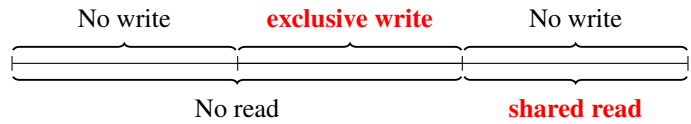
\begin{figure}[htbp]
\begin{center}
\begin{tikzpicture}[]
    \draw[arrows={Bar[]-Bar[]}] (0.0, 0.0) -- (9.0, 0.0);
    \draw[arrows={Bar[]-Bar[]}] (3.0, 0.0) -- (6.0, 0.0);

    \draw[decorate, decoration={brace, raise=5pt}, line width=0.25mm] (0.02, 0.0) -- (2.98, 0.0);
    \node[draw=none] at (1.5, 0.6) {No write};
    \draw[decorate, decoration={brace, raise=5pt}, line width=0.25mm] (3.02, 0.0) -- (5.98, 0.0);
    \node[draw=none] at (4.5, 0.6) {\color{red} \textbf{exclusive write}};
    \draw[decorate, decoration={brace, raise=5pt}, line width=0.25mm] (6.02, 0.0) -- (8.98, 0.0);
    \node[draw=none] at (7.5, 0.6) {No write};

    \draw[decorate, decoration={brace, raise=5pt}, line width=0.25mm] (5.98, 0.0) -- (0.02, 0.0);
    \node[draw=none] at (3.0, -0.6) {No read};
    \draw[decorate, decoration={brace, raise=5pt}, line width=0.25mm] (8.98, 0.0) -- (6.02, 0.0);
    \node[draw=none] at (7.5, -0.6) {\color{red} \textbf{shared read}};
\end{tikzpicture}
\end{center}
\caption{Bucket read/write activity in stages not processed (left), being processed (middle), and processed (right).}
\label{fig:bucket-read-write}
\end{figure}

The dependency discipline makes bucket access conflict-free (\autoref{fig:bucket-read-write}). A bucket not yet released is untouched; a released bucket
is written only by the single worker processing it, which writes only its own storage; and a processed bucket is thereafter read-only, so any number of
workers may read it. Since no bucket is ever read before it is processed, workers never need to coordinate bucket access. The only state that genuinely
requires synchronisation is shared across the whole search -- the upper bound used for pruning and the pool of generated paths -- each guarded by a
single lightweight lock. The algorithm terminates when no buckets remain.

\autoref{alg:pull} summarises the mono-directional algorithm. After computing the step sizes and the dependency counts, labels are
seeded in the source's buckets and every bucket with no remaining dependencies is placed in the ready queue. Workers repeatedly
take a released bucket, pull labels into it from its predecessor buckets (running the dominance check of \autoref{sec:vectorised} as labels are
created), and then call \textsc{RemoveDependency} to release any bucket whose last dependency has just completed. Since each released bucket is
written by exactly one worker and read only after it is processed (\autoref{fig:bucket-read-write}), workers never block one another on bucket data.

\begin{algorithm}[H]
\caption{Parallel mono-directional pull labelling.}
\label{alg:pull}
\begin{algorithmic}
\Procedure{PullLabelling}{}
  \State compute step sizes $\Delta_v$ and dependency counts $DependencyCount[\cdot]$ for all buckets
  \State seed labels in the source's buckets; $\;Q \gets \{\,$buckets with $DependencyCount=0\,\}$
  \While{$Q \neq \emptyset$} \Comment{processed in parallel, one worker per bucket}
     \State $b \gets Q.\Call{pop}{}$
     \State \Call{Extend}{$b$} \Comment{pull feasible labels from predecessor buckets; dominate and store}
     \State \Call{RemoveDependency}{$b$} \Comment{releases newly ready buckets into $Q$}
  \EndWhile
\EndProcedure
\end{algorithmic}
\end{algorithm}

The bi-directional version of \autoref{sec:bidirectional} generalises this driver: it additionally seeds backward labels at the target, and the ready
queue then holds forward, backward, and splice tasks. These are dispatched in a dynamically balanced order that realises the dynamic midpoint, so
\textsc{Extend} above is replaced by the appropriate forward extension, backward extension, or \textsc{Splice} for the task at hand.

\section{Bi-directional Search}
\label{sec:bidirectional}

The number of labels produced is a key performance metric, as it drives both running time and memory, and for a mono-directional search it is typically
exponential in the path length (\autoref{fig:label-counts}(a)). A bi-directional search produces far fewer, because each direction need only reach
``half'' the path length, but labels from the two directions must then be combined at a middle point (\autoref{fig:label-counts}(b), with forward labels
red-shaded and backward labels blue-shaded).

\begin{figure}[htbp]
    \begin{center}
    \begin{minipage}{.48\textwidth}\begin{center}
    \begin{tikzpicture}
    \begin{axis}[grid=none, width=6cm, height=4.5cm, clip=false,
        xmax=12,ymax=8,
        yticklabel=\empty,
        xticklabel=\empty,
        tick style={draw=none},
        axis line style={draw=none},
        enlargelimits]
    \node[red, draw=none] at (74.0, 640.0) {\small Forward};
    \node[draw=none] at (0.0, 850) {Labels};
    \draw[->, line width=0.25mm] (0.0, 0.0) -- (0.0, 800);
    \node[draw=none] at (132, 0.0) {Length};
    \draw[->, line width=0.25mm] (0.0, 0.0) -- (105, 0.0);
    \addplot[red, name path=ff, domain=0:10] {0.5*pow(1.71, x-5)};
    \path[name path=axis] (axis cs:0.0, 0.0) -- (axis cs:10.0, 0.0);
    \addplot[thick, color=red, fill=red, fill opacity=0.05] fill between[of=ff and axis, soft clip={domain=0:10}];
    \end{axis}
    \end{tikzpicture}\\ (a) Mono-directional
    \end{center}\end{minipage}%
    \begin{minipage}{.48\textwidth}\begin{center}
    \begin{tikzpicture}
    \begin{axis}[grid=none, width=6cm, height=4.5cm, clip=false,
        xmax=12,ymax=8,
        yticklabel=\empty,
        xticklabel=\empty,
        axis line style={draw=none},
        tick style={draw=none},
        enlargelimits]
    \node[red, draw=none] at (74.0, 640.0) {\small Forward};
    \node[blue, draw=none] at (28.0, 640.0) {\small Backward};
    \node[draw=none] at (0.0, 850) {Labels};
    \draw[->, line width=0.25mm] (0.0, 0.0) -- (0.0, 800);
    \node[draw=none] at (132, 0.0) {Length};
    \draw[->, line width=0.25mm] (0.0, 0.0) -- (105, 0.0);
    \addplot[red, name path=ff, domain=0:10] {0.5*pow(1.71, x-5)};
    \addplot[blue, name path=fb, domain=-0:10] {0.5*pow(1.71, -x+5)};
    \path[name path=axis] (axis cs:0.0, 0.0) -- (axis cs:10.0, 0.0);
    \addplot[thick, color=blue, fill=blue, fill opacity=0.10] fill between[of=fb and axis, soft clip={domain=5:10}];
    \addplot[thick, color=red, fill=red, fill opacity=0.10] fill between[of=ff and axis, soft clip={domain=0:5}];
    \node[circle, fill=black, inner sep=1.1pt] at (axis cs:5, 0.5) {};
    \node[draw=none] at (axis cs:5, 2.5) {Middle};
    \draw[->, line width=0.25mm] (axis cs:5, 1.9) -- (axis cs:5, 0.75);
    \end{axis}
    \end{tikzpicture}\\ (b) Bi-directional
    \end{center}\end{minipage}
    \end{center}
    \caption{Number of labels as a function of path length. (a) Mono-directional search: the count grows exponentially in the path length. (b)
    Bi-directional search: forward (red) and backward (blue) labels meet at a middle point, roughly halving the length on each side.}
    \label{fig:label-counts}
    \end{figure}

\subsection{Modifications to the Dependency Graph}
\label{sec:splice}

To facilitate a bi-directional search the dependency graph is first modified to include both forward and backward dependencies.
\autoref{fig:dependency-graph-bi}(a) sketches the forward bucket dependencies, which are similar to the mono-directional case, and
\autoref{fig:dependency-graph-bi}(b) the backward case, which is essentially a reversed version of the forward case. To avoid overlapping labels at a
vertex either a forward or a backward bucket at a certain level can be processed -- never both.

\begin{figure}[htbp]
    \begin{center}
    \colorlet{shadecolor}{gray!34}
    \begin{minipage}{.48\textwidth}\begin{center}
    \begin{tikzpicture}[
        nodeWhite/.style={rectangle, draw=, fill=white, very thick, minimum width=8mm, minimum height=6mm},
        nodeGrey/.style={rectangle, draw=, fill=shadecolor, very thick, minimum width=8mm, minimum height=6mm},
        ]
        \node[nodeWhite](f3_1) at (0.0, 0.4){\large $f_3$};
        \node[nodeWhite](f2_1) at (0.0, -0.4){\large $f_2$};
        \node[nodeWhite](f1_1) at (0.0, -1.2){\large $f_1$};
        \node[nodeWhite](b3_1) at (1.0, 0.4){\large $b_3$};
        \node[nodeWhite](b2_1) at (1.0, -0.4){\large $b_2$};
        \node[nodeWhite](b1_1) at (1.0, -1.2){\large $b_1$};
        \node[nodeWhite](f3_2) at (4.0, 0.4){\large $f_3$};
        \node[nodeWhite](f2_2) at (4.0, -0.4){\large $f_2$};
        \node[nodeWhite](f1_2) at (4.0, -1.2){\large $f_1$};
        \node[nodeWhite](b3_2) at (5.0, 0.4){\large $b_3$};
        \node[nodeWhite](b2_2) at (5.0, -0.4){\large $b_2$};
        \node[nodeWhite](b1_2) at (5.0, -1.2){\large $b_1$};
        \draw[->, line width=0.25mm] (0.3, -0.35) -- (3.75, 0.4);
        \draw[->, line width=0.25mm] (0.3, -1.15) -- (3.75, -0.4);
        \draw[->, line width=0.25mm] (0.25, -0.25) -- (0.25, 0.25);
        \draw[->, line width=0.25mm] (0.25, -1.05) -- (0.25, -0.55);
        \draw[->, line width=0.25mm] (4.25, -0.25) -- (4.25, 0.25);
        \draw[->, line width=0.25mm] (4.25, -1.05) -- (4.25, -0.55);
        \node[draw=none] at (1.0, -2.0) {\scriptsize Vertex $x$};
        \node[draw=none] at (5.0, -2.0) {\scriptsize Vertex $y$};
    \end{tikzpicture}\\ (a) Forward dependencies
    \end{center}\end{minipage}%
    \begin{minipage}{.48\textwidth}\begin{center}
    \begin{tikzpicture}[
        nodeWhite/.style={rectangle, draw=, fill=white, very thick, minimum width=8mm, minimum height=6mm},
        nodeGrey/.style={rectangle, draw=, fill=shadecolor, very thick, minimum width=8mm, minimum height=6mm},
        ]
        \node[nodeWhite](f3_1) at (0.0, 0.4){\large $f_3$};
        \node[nodeWhite](f2_1) at (0.0, -0.4){\large $f_2$};
        \node[nodeWhite](f1_1) at (0.0, -1.2){\large $f_1$};
        \node[nodeWhite](b3_1) at (1.0, 0.4){\large $b_3$};
        \node[nodeWhite](b2_1) at (1.0, -0.4){\large $b_2$};
        \node[nodeWhite](b1_1) at (1.0, -1.2){\large $b_1$};
        \node[nodeWhite](f3_2) at (4.0, 0.4){\large $f_3$};
        \node[nodeWhite](f2_2) at (4.0, -0.4){\large $f_2$};
        \node[nodeWhite](f1_2) at (4.0, -1.2){\large $f_1$};
        \node[nodeWhite](b3_2) at (5.0, 0.4){\large $b_3$};
        \node[nodeWhite](b2_2) at (5.0, -0.4){\large $b_2$};
        \node[nodeWhite](b1_2) at (5.0, -1.2){\large $b_1$};
        \draw[<-, line width=0.25mm] (1.3, -0.35) -- (4.75, 0.4);
        \draw[<-, line width=0.25mm] (1.3, -1.15) -- (4.75, -0.4);
        \draw[<-, line width=0.25mm] (1.25, -0.25) -- (1.25, 0.25);
        \draw[<-, line width=0.25mm] (1.25, -1.05) -- (1.25, -0.55);
        \draw[<-, line width=0.25mm] (5.25, -0.25) -- (5.25, 0.25);
        \draw[<-, line width=0.25mm] (5.25, -1.05) -- (5.25, -0.55);
        \node[draw=none] at (1.0, -2.0) {\scriptsize Vertex $x$};
        \node[draw=none] at (5.0, -2.0) {\scriptsize Vertex $y$};
    \end{tikzpicture}\\ (b) Backward dependencies
    \end{center}\end{minipage}
    \end{center}
    \caption{Dependency graph for the bi-directional search: (a) forward bucket dependencies (similar to the mono-directional case) and (b) backward
    bucket dependencies (its reverse). Vertices $x$ and $y$ each carry forward ($f$) and backward ($b$) buckets.}
\label{fig:dependency-graph-bi}
\end{figure}

Secondly, the dependency graph is modified to include splice jobs. A splice job is the process of combining labels in backward
 bucket with relevant labels in forward buckets. That is, forward labels from buckets at other vertices that 
 when connected to backward labels in the backward bucket related to the splice job form feasible paths.

\ifijoc Figure~EC.4 in the online supplement shows\else\autoref{fig:dependency-graph-bi-splice}
shows\fi\ the dependency graph for the splice job. A splice job can be initiated when the backward bucket at the
same level has been processed (black edges) and either the splice job or the forward bucket at the 
previous level (blue edges) and the splice job or the forward bucket at adjacent vertices that should be 
connected to backward labels at this level (red edges). 

These dependencies ensure that every connectable pair of labels is generated and no path is produced twice: a forward label at another vertex has
either already been pulled onward into a forward bucket or splice job, or has reached its midpoint and is ready to be spliced here.

\ifijoc\else
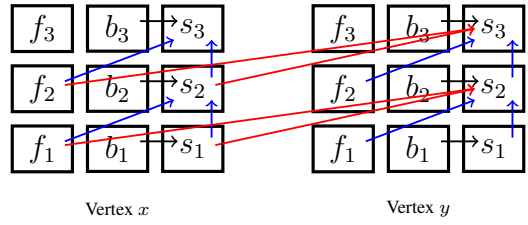
\begin{figure}[htbp]
    \begin{center}
    \colorlet{shadecolor}{gray!34}
    \begin{tikzpicture}[
        nodeWhite/.style={rectangle, draw=, fill=white, very thick, minimum width=8mm, minimum height=6mm},
        nodeGrey/.style={rectangle, draw=, fill=shadecolor, very thick, minimum width=8mm, minimum height=6mm},
        ]
    
        \node[nodeWhite](f3_1) at (0.0, 0.4){\large $f_3$};
        \node[nodeWhite](f2_1) at (0.0, -0.4){\large $f_2$};
        \node[nodeWhite](f1_1) at (0.0, -1.2){\large $f_1$};

        \node[nodeWhite](b3_1) at (1.0, 0.4){\large $b_3$};
        \node[nodeWhite](b2_1) at (1.0, -0.4){\large $b_2$};
        \node[nodeWhite](b1_1) at (1.0, -1.2){\large $b_1$};

        \node[nodeWhite](s3_1) at (2.0, 0.4){\large $s_3$};
        \node[nodeWhite](s2_1) at (2.0, -0.4){\large $s_2$};
        \node[nodeWhite](s1_1) at (2.0, -1.2){\large $s_1$};

        \node[nodeWhite](f3_2) at (4.0, 0.4){\large $f_3$};
        \node[nodeWhite](f2_2) at (4.0, -0.4){\large $f_2$};
        \node[nodeWhite](f1_2) at (4.0, -1.2){\large $f_1$};

        \node[nodeWhite](b3_2) at (5.0, 0.4){\large $b_3$};
        \node[nodeWhite](b2_2) at (5.0, -0.4){\large $b_2$};
        \node[nodeWhite](b1_2) at (5.0, -1.2){\large $b_1$};

        \node[nodeWhite](s3_2) at (6.0, 0.4){\large $s_3$};
        \node[nodeWhite](s2_2) at (6.0, -0.4){\large $s_2$};
        \node[nodeWhite](s1_2) at (6.0, -1.2){\large $s_1$};

        \draw[->, line width=0.25mm] (1.3, 0.5) -- (1.8, 0.5);
        \draw[->, line width=0.25mm] (1.3, -0.3) -- (1.8, -0.3);
        \draw[->, line width=0.25mm] (1.3, -1.1) -- (1.8, -1.1);        
        \draw[->, line width=0.25mm] (5.3, 0.5) -- (5.8, 0.5);
        \draw[->, line width=0.25mm] (5.3, -0.3) -- (5.8, -0.3);
        \draw[->, line width=0.25mm] (5.3, -1.1) -- (5.8, -1.1);

        \draw[->, line width=0.25mm, color=blue] (2.25, -0.25) -- (2.25, 0.25);
        \draw[->, line width=0.25mm, color=blue] (2.25, -1.05) -- (2.25, -0.55);
        \draw[->, line width=0.25mm, color=blue] (6.25, -0.25) -- (6.25, 0.25);
        \draw[->, line width=0.25mm, color=blue] (6.25, -1.05) -- (6.25, -0.55);

        \draw[->, line width=0.25mm, color=blue] (0.3, -0.3) -- (1.75, 0.25);
        \draw[->, line width=0.25mm, color=blue] (0.3, -1.1) -- (1.75, -0.55);
        \draw[->, line width=0.25mm, color=blue] (4.3, -0.3) -- (5.75, 0.25);
        \draw[->, line width=0.25mm, color=blue] (4.3, -1.1) -- (5.75, -0.55);

        \draw[->, line width=0.25mm, color=red] (0.3, -0.35) -- (5.75, 0.4);
        \draw[->, line width=0.25mm, color=red] (0.3, -1.15) -- (5.75, -0.4);

        \draw[->, line width=0.25mm, color=red] (2.3, -0.35) -- (5.75, 0.4);
        \draw[->, line width=0.25mm, color=red] (2.3, -1.15) -- (5.75, -0.4);        
        
        \node[draw=none] at (1.0, -2.0) {\scriptsize Vertex $x$};
        \node[draw=none] at (5.0, -2.0) {\scriptsize Vertex $y$};
    \end{tikzpicture}
    \end{center}
    \caption{Dependency graph for bi-directional search for splice job dependencies.}
\label{fig:dependency-graph-bi-splice}
\end{figure}
\fi

The release rules for the forward bucket, backward bucket, and splice job at vertex $y$, level $n$, with connected vertex $x$ are summarised below.
They are written for the schematic equal-bucket case; for buckets of differing sizes the indices $n\pm1$ are replaced by the predecessor/successor
bucket maps, as noted after the display:

\[
\begin{aligned}
\Call{Release}{f^y_n} &\gets \Call{processed}{f^y_{n-1}} \land \Call{processed}{f^x_{n-1}} \\
\Call{Release}{b^y_n} &\gets \Call{processed}{b^y_{n+1}} \land \Call{processed}{b^x_{n+1}} \\
\Call{Release}{s^y_n} &\gets \Call{processed}{b^y_n} \land \left(\Call{processed}{f^y_{n-1}} \lor \Call{processed}{s^y_{n-1}}\right) \\
&\phantom{{}\gets{}} \land \left(\Call{processed}{f^x_{n-1}} \lor \Call{processed}{s^x_{n-1}}\right)
\end{aligned}
\]

The conjunction over the connected vertex $x$ runs over all in-edges $(x,y)$; the backward rule mirrors $\pi$ towards larger resource values, and the
splice rule combines the forward predecessor map with the same-vertex backward bucket. The splice job itself is given in \autoref{alg:splice}: for each
backward label it joins, across every in-edge, exactly the forward labels whose resource still leaves room for the edge, exploiting that forward labels
are stored in increasing $\rho$ to stop early. Here $\rho(b)$ is the backward label's $\rho$ on the forward scale (its value at the splice vertex), so
$f$ joins $b$ across edge $e$ only if $\rho(f)+q^\rho_e\le\rho(b)$; $\mathit{rc}(e)$ is the reduced cost of $e$, $\mathit{UB}$ the current pruning upper
bound, and $\Call{FeasibleSplice}{f,b,e}$ tests joint resource feasibility of $f\oplus e\oplus b$ (including $ng$-memory compatibility).

\begin{algorithm}[H]
\caption{Splice: concatenate forward and backward labels at a vertex.}
\label{alg:splice}
\begin{algorithmic}
\Procedure{Splice}{vertex $y$, backward bucket $k$}
  \ForAll{backward labels $b$ in bucket $k$ of $y$}
     \ForAll{feasible in-edges $e=(x,y)$ with consumption $t = q^\rho_{e}$}
        \ForAll{forward labels $f$ of $x$ in increasing $\rho$}
           \If{$\rho(f) + t > \rho(b)$} \State \textbf{break} \Comment{no later $f$ can join $b$ across $e$} \EndIf
           \State $c \gets cost(f) + cost(b) + \mathit{rc}(e)$
           \If{$\Call{FeasibleSplice}{f,b,e} \land c < \mathit{UB}$}
              \State \Call{StorePath}{$c$, $f \,{\oplus}\, e \,{\oplus}\, b$}
           \EndIf
        \EndFor
     \EndFor
  \EndFor
\EndProcedure
\end{algorithmic}
\end{algorithm}

We now argue that this enumerates every feasible path exactly once.

\begin{lemma}[Splice completeness]
\label{lem:splice}
Assume the release rules above. Then every feasible source--target path is produced exactly once: either by direct extension to the target, or by
exactly one splice of a forward and a backward label.
\end{lemma}
\begin{proof}[Proof sketch]
This is the bi-directional half-way principle \citep{salani2006,tilk2017} realised through the dependency graph. Let $P=(v_0,\dots,v_m)$ be a
feasible source--target path and let $\rho_i$ be the accumulation of the monotone resource $\rho$ at $v_i$, measured on the forward scale. Since
$q^\rho_e>0$ on every edge, $\rho_0<\rho_1<\dots<\rho_m$, and the forward partial $L^f_i$ and the backward partial $L^b_i$ of $P$ at $v_i$ both
correspond to the single bucket of $v_i$ that contains $\rho_i$: the resource consumed from the source to $v_i$ is the same value in either direction,
and we read $\rho(b)$ on this forward scale throughout. Three facts then combine.
\emph{(1) Prefix/suffix.} A forward partial $L^f_i$ exists only if its parent $L^f_{i-1}$ does, because forward extension creates $L^f_i$ by pulling
from $L^f_{i-1}$ across $(v_{i-1},v_i)$; hence the indices for which $L^f_i$ is created form a prefix $\{0,\dots,p\}$, and symmetrically the backward
partials form a suffix $\{q,\dots,m\}$.
\emph{(2) Coverage in exactly one direction.} Within a vertex the release rules of \autoref{sec:splice} process forward buckets in increasing $\rho$ and
backward buckets in decreasing $\rho$, so a vertex's forward-covered levels are down-closed, its backward-covered levels are up-closed, and the
never-both invariant keeps the two disjoint. This coverage is moreover consistent \emph{along} $P$: if $v_i$ is backward-covered then so is $v_{i+1}$,
because the backward partial at $v_i$ is created by extending the backward partial at $v_{i+1}$ across $(v_i,v_{i+1})$, so $v_{i+1}$ must be
backward-covered before $v_i$ can be; equivalently, forward coverage is inherited from $v_{i+1}$ back to $v_i$. Hence the covered directions along $P$
form a block of forward-covered indices followed by a block of backward-covered indices with a single transition, and -- since each level is claimed by
whichever search reaches it first and the two searches together exhaust the levels -- the two blocks tile $\{0,\dots,m\}$ with no gap and no overlap.
Their boundary is the dynamic half-way cut of \citet{salani2006,tilk2017}, emerging from the bucket processing order rather than fixed a priori. Because
$\rho$ is strictly monotone along $P$, every vertex $v_i$ lies on one definite side of this cut, so each index $i\in\{0,\dots,m\}$ is forward-covered or
backward-covered, but not both.
\emph{(3) Adjacency.} Thus $\{0,\dots,p\}$ and $\{q,\dots,m\}$ are disjoint and together cover $\{0,\dots,m\}$, which forces $q=p+1$. Hence $P$ splits
at the unique edge $e=(v_p,v_{p+1})$ with $f=L^f_p$ and $b=L^b_{p+1}$, and the join is feasible because $\rho(f)+q^\rho_e=\rho_{p+1}=\rho(b)$. (If the
prefix already reaches the target, $p=m$ and $P$ is produced directly by forward extension to the target instead.)
Finally $b$ lies in exactly one backward bucket $(v_{p+1},k)$, the argument of exactly one splice job $\textsc{Splice}(v_{p+1},k)$ released once by its
dependency counter, and within it the loop enumerates the forward labels of $v_p$ with $\rho(f)+q^\rho_e\le\rho(b)$ once. The split edge is unique, so
$P$ is not double counted; and step~(2) guarantees that $f$ and $b$ both exist when the job runs, so $P$ is not missed. Hence $P$ is produced exactly once.
\end{proof}

To illustrate how the forward and backward buckets and splice jobs are processed within a 
vertex, \ifijoc Figure~EC.1 in the online supplement shows\else\autoref{fig:bucket-update} shows\fi\ the bucket update activity during its lifetime.
Red buckets are not ready for processing, blue buckets are ready for processing and greyed out 
buckets are processed. This assumes that the release rules are used to determine when a bucket is 
ready for processing including dependencies from adjacent vertices. Forward
and backward buckets are processed in parallel from bottom to top and vice versa. Splice jobs are 
initiated when the backward bucket at the same level has been processed and goes upward from the 
level of the last backward bucket. The remaining jobs in red are never initiated or processed.

\ifijoc\else
\input{fig_bucket_update}
\fi

\subsection{Dynamic Midpoint}

To minimise the total number of labels, the midpoint should balance the number of forward and backward labels. A fixed midpoint stops extending a
label once its monotone resource $\rho$ reaches a preset value, e.g.\ half the resource bound $\lceil b^\rho/2 \rceil$.

Following \citet{tilk2017}, we instead use a dynamic midpoint. Keeping separate queues for forward, backward, and splice jobs and counting the forward
and backward buckets processed, the scheduler always dispatches a ready splice job first and otherwise the direction with the smaller processed count.
No midpoint value is ever fixed: each direction advances until the two have processed equally many buckets, so the balance emerges from the processing
order.

\section{Vectorised Resource Comparison}
\label{sec:vectorised}

The dominance check -- comparing a new label against the stored labels of its vertex to decide whether it is dominated -- is usually done one label at a
time; it can instead be vectorised to compare against many labels at once.

\subsection{Dominance}

For a label $L$ we write $v(L)$ for its vertex, $cost(L)$ for its accumulated cost, and $r(L)$ for its accumulation of resource $r \in R$.
Let $\mathcal{E}(L)$ be the set of valid extensions of $L$, i.e.\ the paths from $v(L)$ to the target vertex that extend $L$ feasibly, and let
$cost(L+\epsilon)$ denote the cost of $L$ extended by $\epsilon \in \mathcal{E}(L)$. The dominance criterion can be defined as follows:

\begin{definition}[Dominance]
\label{def:dominance}
Let $L_1$ and $L_2$ be two labels at the same vertex, $v(L_1)=v(L_2)$. Label $L_1$ dominates label $L_2$ if every feasible extension of $L_2$ is also
a feasible extension of $L_1$ at no greater cost:
$$\mathcal{E}(L_2) \subseteq \mathcal{E}(L_1) \quad\text{and}\quad cost(L_1 + \epsilon) \leq cost(L_2 + \epsilon) \ \ \forall \epsilon \in \mathcal{E}(L_2).$$
\end{definition}

When $L_1$ dominates $L_2$, every path through $L_2$ has a no-costlier counterpart through $L_1$, so $L_2$ can be discarded. As not all extensions of
a label are known, a sufficient condition is used to determine dominance.

\begin{lemma}[Sufficient condition]
\label{lemmaSufficientCondition}
Assume each REF is monotone non-decreasing in the incoming resource state (so that resets, which are not, are excluded). Then label $L_1$ dominates label $L_2$ if:
\begin{center}
$
\begin{aligned}
    v(L_1) & = v(L_2) \\
    cost(L_1) & \leq cost(L_2) \\
    r(L_1) & \leq r(L_2) : \forall r \in R
\end{aligned}
$
\end{center}
\end{lemma}

Based on the sufficient condition a natural total order of labels can be defined. Let $\rho$ be the monotone resource on which buckets are built
and let the remaining resources be $r_2,\dots,r_{|R|}$. For two labels at the same vertex we order them lexicographically by
$$L_1 \leq^{lex} L_2 \iff \big(\rho(L_1),\, cost(L_1),\, r_2(L_1),\dots\big) \ \text{precedes or equals}\ \big(\rho(L_2),\, cost(L_2),\, r_2(L_2),\dots\big),$$
i.e.\ first by the monotone resource, then by cost, then by the remaining resources in a fixed order. On distinct labels this is a total order:
otherwise incomparable labels are ordered deterministically by the first coordinate on which they differ, and two labels that agree on every coordinate
are duplicates (mutually dominating), one of which is discarded on creation.

\begin{lemma}[Dominance order]
\label{lemmaDominanceOrder}
Under the assumptions of Lemma~\ref{lemmaSufficientCondition}, the lexicographic order $\leq^{lex}$ is consistent with dominance, where
$\preceq_{dom}$ denotes the sufficient-condition order of \autoref{lemmaSufficientCondition}:
$$L_1 \preceq_{dom} L_2 \Rightarrow L_1 \leq^{lex} L_2.$$
\end{lemma}

That is, if $L_2$ is dominated by $L_1$ then $L_1$ comes first in the order. We now make precise the sense in which the pull order keeps the
stored labels a Pareto frontier, which is what justifies immutable storage.

\begin{lemma}[Pareto frontier by construction]
\label{lem:immutable}
Fix a vertex and suppose all REFs are non-decreasing and the monotone resource $\rho$ strictly increases along every edge
(\autoref{sec:preliminaries}). Suppose its labels are created in non-decreasing $\leq^{lex}$ order, which the pull schedule guarantees: the
vertex's forward buckets are processed in increasing $\rho$, and the batch pulled into one bucket is inserted in $\leq^{lex}$ order (this per-batch sort
is required: without it a newly created label could be dominated by another label of the same batch that is not yet stored, and testing only against the
already-stored labels would be unsound). Then a label
that is not dominated when created is never dominated by a label created later. Consequently it suffices to test each new label against the
already-stored labels of the vertex, and a stored label is never removed.
\end{lemma}
\begin{proof}[Proof sketch]
Let $L$ be stored and let $L'$ be created later at the same vertex, so $L \leq^{lex} L'$ and in particular $\rho(L)\le\rho(L')$. If
$\rho(L') > \rho(L)$ then $L'$ violates $r(L')\le r(L)$ on the resource $\rho$, so by \autoref{lemmaSufficientCondition} $L'$ cannot dominate $L$.
If $\rho(L')=\rho(L)$ and $L'\neq L$, then $L \leq^{lex} L'$ means that on the first coordinate among $(cost,r_2,\dots)$ where they differ
$L' > L$, so $L'$ again exceeds $L$ on cost or on some resource and cannot dominate $L$. Finally, if $L'$ equals $L$ in every coordinate then $L'$
is dominated by the stored $L$ and is discarded on creation. In all cases no later label dominates a stored label. Since therefore every potential
dominator of a label is created no later than the label itself, comparing a new label only against the already-stored labels is sufficient, and an
accepted label is permanent. The backward search is symmetric, ordering by decreasing $\rho$ and ascending cost.
\end{proof}

This directly addresses the concern that the order between mutually incomparable labels is arbitrary: $\leq^{lex}$ is a fixed total order, and
incomparable labels are simply both retained. Because a stored label is never removed (only an incoming label can be discarded on creation),
labels live at a fixed position in memory and can be stored as immutable objects, avoiding any moving or copying of labels.

\subsection{Structure of Resources}

The controlled creation order lets us store labels as a \emph{structure of arrays} rather than an \emph{array of structures} (\autoref{fig:label-layout}).
Writing a label's members (its resources, e.g.\ cost, time, load) as $\{a,b,c\}$: the array-of-structures layout keeps each label's members together and
lays labels out consecutively, whereas the structure-of-arrays layout stores each member in its own contiguous array across all labels. The latter is
exactly what the dominance check needs, since it compares one member across many labels at a time -- a contiguous scan with few cache misses -- and is
what the vectorised test consumes.

\begin{figure}[htbp]
\begin{center}
\begin{tikzpicture}[
    node/.style={rectangle, draw=, fill=white, very thick, minimum width=10mm, minimum height=10mm},
    ]
    \node[node](a0) at (0.0, 0.0){\LARGE $a_0$};
    \node[node](b0) at (1.0, 0.0){\LARGE $b_0$};
    \node[node](c0) at (2.0, 0.0){\LARGE $c_0$};
    \node[node](a1) at (3.0, 0.0){\LARGE $a_1$};
    \node[node](b1) at (4.0, 0.0){\LARGE $b_1$};
    \node[node](c1) at (5.0, 0.0){\LARGE $c_1$};
    \node[node](a2) at (6.0, 0.0){\LARGE $a_2$};
    \node[node](b2) at (7.0, 0.0){\LARGE $b_2$};
    \node[node](c2) at (8.0, 0.0){\LARGE $c_2$};
    \node[draw=none] at (9.0, 0.0) {...};
    \node[draw=none] at (4.0, -1.0) {\small (a) Array of structures};
\end{tikzpicture}

\vspace{1.2em}

\begin{tikzpicture}[
    node/.style={rectangle, draw=, fill=white, very thick, minimum width=10mm, minimum height=10mm},
    ]
    \node[node](a1) at (0.0, 0.1){\LARGE $a_0$};
    \node[node](a2) at (1.0, 0.1){\LARGE $a_1$};
    \node[node](a3) at (2.0, 0.1){\LARGE $a_2$};
    \node[draw=none] at (3.0, 0.1) {...};
    \node[node](b1) at (0.0, -1.0){\LARGE $b_0$};
    \node[node](b2) at (1.0, -1.0){\LARGE $b_1$};
    \node[node](b3) at (2.0, -1.0){\LARGE $b_2$};
    \node[draw=none] at (3.0, -1.0) {...};
    \node[node](c1) at (0.0, -2.1){\LARGE $c_0$};
    \node[node](c2) at (1.0, -2.1){\LARGE $c_1$};
    \node[node](c3) at (2.0, -2.1){\LARGE $c_2$};
    \node[draw=none] at (3.0, -2.1) {...};
    \node[draw=none] at (1.0, -3.1) {\small (b) Struct of arrays};
\end{tikzpicture}
\end{center}
\caption{Label memory layout: (a) array of structures, where each label's members are contiguous, versus (b) struct of arrays, where each member's
values are contiguous across labels -- the layout the vectorised dominance test consumes.}
\label{fig:label-layout}
\end{figure}

\subsection{Vector Instruction Set}

\begin{figure}[htbp]
\begin{center}
\begin{tikzpicture}[
    node/.style={rectangle, draw=, fill=white, very thick, minimum width=10mm, minimum height=10mm},
    ]
    \node[node](a0) at (0.0, 0.0){\LARGE $a_0$};
    \node[node](a1) at (1.0, 0.0){\LARGE $a_1$};
    \node[node](a2) at (2.0, 0.0){\LARGE $a_2$};
    \node[draw=none] at (3.0, 0.0) {...};
    \node[node](a31) at (4.0, 0.0){\LARGE $a_{W-1}$};
    \node[draw=none] at (5.4, 0.0) {\huge $\leq$};
    \node[node](an0) at (6.6, 0.0){\LARGE $a^\star$};
    \node[node](an1) at (7.6, 0.0){\LARGE $a^\star$};
    \node[node](an2) at (8.6, 0.0){\LARGE $a^\star$};
    \node[draw=none] at (9.6, 0.0) {...};
    \node[node](an31) at (10.6, 0.0){\LARGE $a^\star$};

    \node[draw=none] at (5.4, -0.8) {\huge $\land$};

    \node[node](bn0) at (0.0, -1.6){\LARGE $b_0$};
    \node[node](bn1) at (1.0, -1.6){\LARGE $b_1$};
    \node[node](bn2) at (2.0, -1.6){\LARGE $b_2$};
    \node[draw=none] at (3.0, -1.6) {...};
    \node[node](bn31) at (4.0, -1.6){\LARGE $b_{W-1}$};
    \node[draw=none] at (5.4, -1.6) {\huge $\leq$};
    \node[node](bn0b) at (6.6, -1.6){\LARGE $b^\star$};
    \node[node](bn1b) at (7.6, -1.6){\LARGE $b^\star$};
    \node[node](bn2b) at (8.6, -1.6){\LARGE $b^\star$};
    \node[draw=none] at (9.6, -1.6) {...};
    \node[node](bn31b) at (10.6, -1.6){\LARGE $b^\star$};

    \node[draw=none] at (5.4, -2.4) {\huge $\land$};

    \node[node](cn0) at (0.0, -3.2){\LARGE $c_0$};
    \node[node](cn1) at (1.0, -3.2){\LARGE $c_1$};
    \node[node](cn2) at (2.0, -3.2){\LARGE $c_2$};
    \node[draw=none] at (3.0, -3.2) {...};
    \node[node](cn31) at (4.0, -3.2){\LARGE $c_{W-1}$};
    \node[draw=none] at (5.4, -3.2) {\huge $\leq$};
    \node[node](cn0c) at (6.6, -3.2){\LARGE $c^\star$};
    \node[node](cn1c) at (7.6, -3.2){\LARGE $c^\star$};
    \node[node](cn2c) at (8.6, -3.2){\LARGE $c^\star$};
    \node[draw=none] at (9.6, -3.2) {...};
    \node[node](cn31c) at (10.6, -3.2){\LARGE $c^\star$};
\end{tikzpicture}
\end{center}
\caption{Vectorised dominance check.}
\label{fig:vectorised-dominance}
\end{figure}
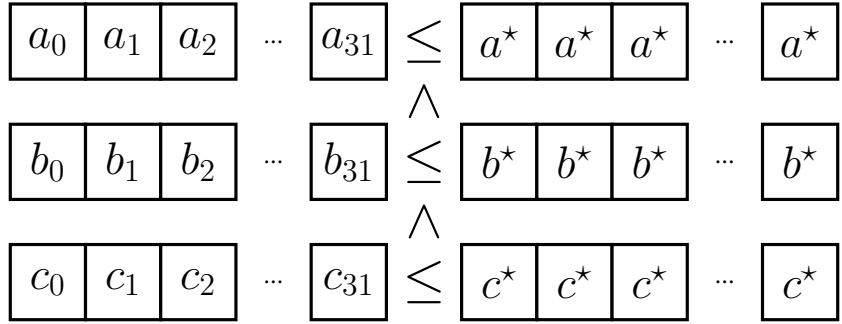

Using vector instructions, label properties can be compared in parallel.
The resource values are $32$-bit integers, so a $256$-bit AVX2 register, used throughout our experiments, holds $W=8$ of them (wider registers such as AVX-512 would hold proportionally more).
The candidate label's value for a resource is broadcast across a register and compared in parallel with the same resource of $W$ other labels.
\autoref{fig:vectorised-dominance} shows a vectorised dominance check.

\autoref{alg:dominance} makes the test precise. The stored labels of a vertex are laid out as a structure of arrays, one contiguous array per
resource, in the lexicographic order $(\rho, cost, \dots)$ of \autoref{lemmaDominanceOrder}. The candidate $L^\star$ is tested against the stored
labels in blocks of $W$ lanes. For a block, the per-lane mask starts all-true and is intersected ($\land$) across resources by one vector comparison
per resource, so a lane survives iff that stored label is $\le L^\star$ on every resource. The cost is held as a scalar ($64$-bit, not packed into the
vectors); when a block reports a surviving lane, a short scalar pass confirms domination by checking that some surviving lane also has
$cost \le cost(L^\star)$. The two SIMD reductions are thus a logical \emph{and} across the resource dimensions and a logical \emph{or} across the lanes
of a block, with the scalar cost giving the final comparison. The lexicographic ordering additionally lets the scan terminate once no remaining label
can beat $L^\star$.

\begin{algorithm}[H]
\caption{Vectorised dominance test.}
\label{alg:dominance}
\begin{algorithmic}
\Procedure{IsDominated}{candidate $L^\star$, stored labels $\mathcal{L}$ (struct-of-arrays)}
  \ForAll{blocks of $W$ lanes in $\mathcal{L}$}
     \If{$\rho(\text{first lane of block}) > \rho(L^\star)$} \State \textbf{break} \Comment{labels stored in increasing $\rho$; no later label can dominate $L^\star$} \EndIf
     \State $\mathit{mask} \gets \mathbf{true}^{W}$
     \ForAll{resources $r$}
        \State $\mathit{mask} \gets \mathit{mask} \ \land\ \big(\,\mathit{block}_r \le \mathrm{broadcast}(r(L^\star))\,\big)$ \Comment{one SIMD compare per resource}
     \EndFor
     \If{$\Call{any}{\mathit{mask}}$} \Comment{a stored label beats $L^\star$ on all resources}
        \ForAll{lanes $j$ with $\mathit{mask}_j$ set}
           \If{$cost(\text{lane } j) \le cost(L^\star)$} \State \Return \textbf{true} \Comment{that label dominates $L^\star$} \EndIf
        \EndFor
     \EndIf
  \EndFor
  \State \Return \textbf{false}
\EndProcedure
\end{algorithmic}
\end{algorithm}

For binary visit-set resources the per-resource comparison is the packed bitwise test $(\,\mathit{block}_r \mathbin{\&} {\sim} r(L^\star) = 0\,)$
rather than $\le$, evaluated in the same vectorised manner. Thus the test reports $L^\star$ dominated iff $\exists L \in \mathcal{L}: L \preceq_{dom} L^\star$, but
the test is performed $W$ labels at a time.

\section{Experiments}
\label{sec:experiments}

To evaluate the performance of the proposed algorithmic optimisations the following experimental setup is considered,
i) a vanilla version of the algorithm,
ii) a version with parallel pulling,
iii) a version with bi-directional search with dynamic midpoint,
iv) a version with vectorised dominance comparison, and
v) a version with all optimisations combined. The configurations are given in \autoref{tab:configurations}.
The \emph{vanilla} (Base) version is the mono-directional pull labelling algorithm of \autoref{sec:mono} run on a single thread,
with scalar (non-vectorised) one-by-one dominance checks; it is the common ancestor of all other configurations, so that each row of
\autoref{tab:configurations} isolates the effect of exactly one optimisation. The parallel configurations use all $32$ hardware threads of the
machine (the $16$ cores with two-way SMT enabled; one thread per released bucket, up to the thread count).

\begin{table}[htbp]
\begin{center}
\begin{tabular}{lccc}
\toprule
Configuration & Parallel & Bi-directional & Vectorised \\
\midrule
Base & & & \\
Parallel & \checkmark & & \\
Bi-directional & & \checkmark & \\
Vectorised & & & \checkmark \\
All & \checkmark & \checkmark & \checkmark \\
\bottomrule
\end{tabular}
\vspace{1em}
\caption{Algorithm configurations}
\label{tab:configurations}
\end{center}
\end{table}

All experiments were performed on a machine with an AMD Ryzen 9 3950X 16-Core Processor ($32$ hardware threads with two-way SMT) with AVX2 support running Ubuntu 24.04.
The algorithm is implemented in C++ and is part of the Flowty proprietary solver. Each instance--configuration pair is timed in a single run; because the effects we report are large (from
roughly $2\times$ to over $200\times$) relative to the few-percent run-to-run variation of compute-bound code, and because runtimes are aggregated
over $56$ instances with a shifted geometric mean that damps per-run noise, the reported comparisons are robust to timing variability. The one place
where single-run noise is non-negligible is the sub-second N8 instances, where it inflates ratios; we therefore do not read much into the very large
N8 factor.

\subsection{Instances}
Instances are generated as pricing problems of the root nodes when solving the Solomon instances \citep{solomon1987}
 for the vehicle routing problem with time windows (VRPTW) using column generation.
 
The 56 original instances are divided into six classes: R1, R2, C1, C2, RC1, and RC2
 representing different combinations of random and clustered customer locations, as well as short and long time windows.
For each original instance we create 3 instances based on the neighbourhood size.
 The N$x$ in the instance name refers to the maximum size $x \in \{8, 16, 24\}$ of the neighbourhood.

 Neighbourhoods are grown dynamically: across column generation iterations they are enlarged to remove cycles up to the maximum size, and once a valid
 lower bound is reached they are trimmed to the current solution's path before growing resumes. For each original instance and each maximum size we keep
 the single pricing instance at the $3600$-second root-node cutoff (or the last one solved if column generation converges earlier); small instances may
 not reach the larger maximum sizes. All $56 \times 3 = 168$ pricing instances are pairwise distinct: hashing the canonical content of each (vertices,
 reduced-cost edges, and neighbour lists) gives $168$ distinct hashes, reproducible via \texttt{scripts/dedup\_check.py}.

The bucket grid is a structural property of an instance: with time as the monotone resource, vertex $v$ has step size
$\Delta_v = \min\{\,\text{time}(e) : e \text{ incident to } v\,\}$ and therefore $N_v = \lfloor (b_v - a_v)/\Delta_v\rfloor + 1$ buckets.
\ifijoc Table~EC.4 in the online supplement reports\else\autoref{tab:bucket_stats} reports\fi\ the resulting average number of buckets per vertex per
class; it averages $22.0$ over all instances and, as expected, does not depend on the maximum neighbourhood size (which changes the labels, not the
bucket grid). Classes with long time windows (C2, R2, RC2) yield substantially more buckets per vertex than the short-window classes (R1, RC1). These
are allocated buckets; the number of \emph{non-empty} buckets at run time is smaller and depends on the dual values of the pricing iteration.

\ifijoc\else
\begin{table}[htbp]
\begin{center}
\begin{tabular}{lc}
\toprule
Class & Buckets per vertex \\
\midrule
C1 & 14.8 \\
C2 & 41.2 \\
R1 & 6.9 \\
R2 & 35.5 \\
RC1 & 7.0 \\
RC2 & 30.1 \\
\midrule
All & 22.0 \\
\bottomrule
\end{tabular}
\vspace{1em}
\caption{Average number of buckets per vertex by instance class. The count is structural, i.e.\ the time window divided by the smallest time consumption over the incident edges, $N_v=\lfloor (b_v-a_v)/\Delta_v\rfloor+1$; it is independent of the maximum neighbourhood size, which changes the labels but not the bucket grid.}
\label{tab:bucket_stats}
\end{center}
\end{table}

\fi

Instances are available at \citep{spoorendonk2025}.

\subsection{Comparison with a state-of-the-art implementation}
\label{sec:comparison}

To place the pull algorithm against the state of the art we compare it with the bucket graph labelling algorithm of \citet{sadykov2020},
the labelling engine at the core of the leading exact branch-cut-and-price solvers for vehicle routing. We use \texttt{bgspprc} \citep{bgspprcpaper,bgspprc},
an independent open-source (MIT) implementation of that algorithm. We are upfront that \texttt{bgspprc} is our own implementation rather than a
third party's; to support its use as a faithful baseline it is publicly available, unit tested, validated against the reference optima of the
instances, and cross-validated against the open PathWyse solver \citep{pathwyse}, which it outperforms by $1.3$--$2.35\times$ \citep{bgspprcpaper}.
\texttt{bgspprc} is therefore a strong open baseline rather than a weak one. It is run in its best parallel configuration (parallel, bi-directional,
with SIMD dominance), which is the natural counterpart of the pull algorithm with all optimisations.

Both solvers were run on the same machine (\autoref{sec:experiments}) on the same $56 \times 3$ instances, each in its best parallel configuration.
\texttt{bgspprc}, like push labelling in general, parallelises the two search directions and uses standard parallel routines such as parallel sorting.
Bucket-level parallelism is awkward for push, since extending one bucket writes labels into several successor buckets that other threads also write;
pull avoids this by construction, each bucket being processed by the single thread that reads only finalised predecessors and writes only its own
storage. The advantage is thus this contention-free, cache-local access pattern rather than a higher nominal degree of parallelism -- indeed the acyclic
dependency graph can \emph{limit} how many buckets are ready at once -- and comparing each solver in its best parallel configuration measures exactly
that. Following the companion study, we report the shifted geometric mean (shift $1$s) per neighbourhood size with a $120$s timeout substituted, and
speed-up $(\text{bg}_{\text{sgm}}+1)/(\text{pull}_{\text{sgm}}+1)$ (\autoref{tab:comparison}). The pull algorithm is $1.9\times$--$2.4\times$ faster, the
gap narrowing monotonically as $n_g$ grows ($2.40 \to 2.10 \to 1.92$), suggesting a fixed per-instance advantage rather than one that scales. The
aggregate margin is shift-dependent -- $1.9$--$2.4\times$ at the $1$s shift, narrowing to about $1.3\times$ at a $10$s shift that weights the
budget-capped instances more heavily -- but the instance-level win rate is shift-independent: pull is faster on $164$ of the $168$ instances. The four
exceptions are hard long-time-window instances at N16/N24, where \texttt{bgspprc} is strictly faster on two (C204 at N24, RC204 at N16) and both solvers
hit the budget on the other two. At $n_g=24$ pull leaves three of $56$ instances unsolved within $120$s (several minutes on its worst, C204 and R208)
against \texttt{bgspprc}'s four -- a heavier tail consistent with the acyclic-dependency restriction limiting parallelism when buckets are small.
Reproducible via \texttt{scripts/comparison.py} and the bundled \texttt{scripts/bgspprc\_para\_bidir\_vec.csv}.

\begin{table}[htbp]
\centering
\begin{tabular}{lrrrrrrr}
\toprule
$n_g$ & \multicolumn{3}{c}{Pull (All)} & \multicolumn{3}{c}{bgspprc} & Speed-up \\
 & sgm & mean & solved & sgm & mean & solved & ($\times$) \\
\midrule
8 & 0.203 & 0.406 & 56/56 & 1.884 & 5.640 & 56/56 & 2.40$\times$ \\
16 & 0.526 & 3.010 & 56/56 & 2.200 & 8.783 & 56/56 & 2.10$\times$ \\
24 & 0.873 & 9.123 & 53/56 & 2.588 & 15.433 & 52/56 & 1.92$\times$ \\
\bottomrule
\end{tabular}
\vspace{1em}
\caption{Same-hardware comparison of the pull algorithm (All optimisations) with the open bucket-graph solver \texttt{bgspprc} (best parallel mode) on the $56$ Solomon instances per $n_g$. Times in seconds; sgm is the shifted geometric mean ($\exp(\overline{\log(t+1)})-1$); the $120$s timeout is substituted, never dropped. Speed-up $=(\text{bg}_{\text{sgm}}+1)/(\text{pull}_{\text{sgm}}+1)$, matching the convention of the companion paper.}
\label{tab:comparison}
\end{table}

\subsection{Performance}
Experiments are run for each instance and configuration. \ifijoc Runtime details for every instance are provided in the online
supplement.\else Runtime details are given in \autoref{sec:runtime-details}.\fi

We aggregate runtimes with the shifted geometric mean, a standard benchmarking measure that limits the influence of the fastest and slowest instances:
for runtimes $t_1,\dots,t_n$ and shift $\sigma$ it is $\left(\prod_{i=1}^{n}(t_i+\sigma)\right)^{1/n}-\sigma$. We use $\sigma = 1$ second, matching the
companion study. The tables report this value unscaled (seconds), scaled by the Base configuration, and as its reciprocal speed-up.

The performance profile (\ifijoc Figure~EC.5 of the online supplement\else\autoref{fig:performance-profile}\fi) shows the combination of all
optimisations performing best, never more than a factor $2$ from the best configuration and up to $274\times$ faster than the baseline on the instance
that benefits most (RC208), with parallelisation and bi-directional search the largest contributors.

\ifijoc\else
\begin{figure}[htbp]
    \centering
    \includegraphics[width=0.8\textwidth]{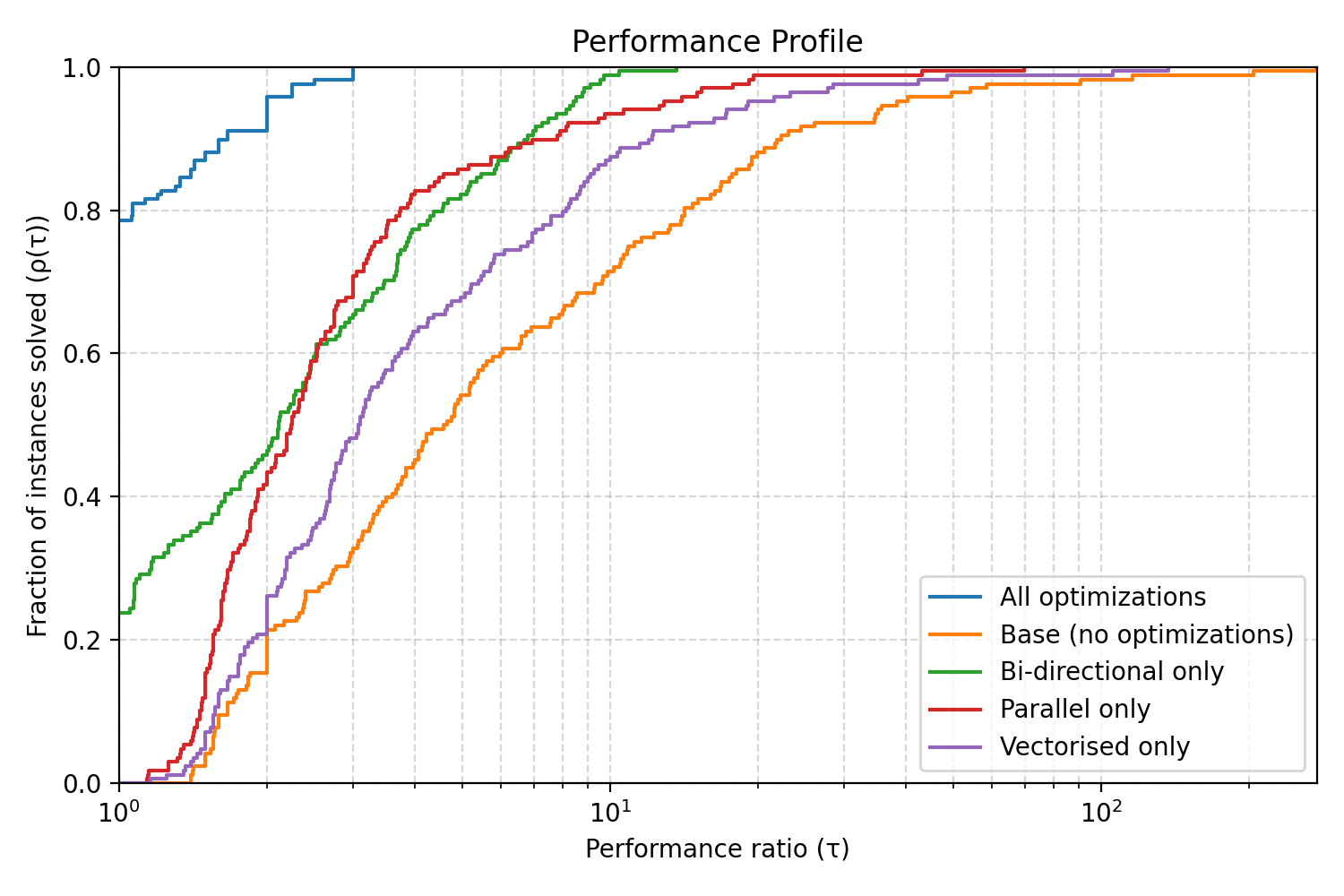}
    \caption{Performance profile of the algorithm with different configurations}
    \label{fig:performance-profile}
  \end{figure}
\fi

The per-neighbourhood-size gains (\autoref{tab:runtime_ng_stats}) are consistent across sizes: around $1.8$--$2.4\times$ from parallelisation,
$1.7\times$ from bi-directional search, $1.4\times$ from vectorised dominance, and around $4\times$ from all combined.

\begin{table}[htbp]
\centering
\begin{tabular}{lrrrrrr}
\toprule
NG & Base & Parallel & Bi-directional & Vectorised & All \\
\midrule
\multicolumn{6}{l}{Runtime (s)} \\
\midrule
8 & 0.957 & 0.403 & 0.533 & 0.628 & 0.203 \\
16 & 2.142 & 1.209 & 1.262 & 1.601 & 0.526 \\
24 & 3.575 & 2.025 & 2.011 & 2.598 & 0.935 \\
\midrule
\multicolumn{6}{l}{Normalized} \\
\midrule
8 & 4.722 & 1.990 & 2.627 & 3.099 & 1.000 \\
16 & 4.069 & 2.297 & 2.397 & 3.042 & 1.000 \\
24 & 3.822 & 2.165 & 2.149 & 2.778 & 1.000 \\
\midrule
\multicolumn{6}{l}{Speedup} \\
\midrule
8 & 1.000 & 2.373 & 1.797 & 1.524 & 4.722 \\
16 & 1.000 & 1.772 & 1.698 & 1.338 & 4.069 \\
24 & 1.000 & 1.766 & 1.778 & 1.376 & 3.822 \\
\bottomrule
\end{tabular}
\vspace{1em}
\caption{Runtime shifted (1s) geometric mean (s) unscaled, scaled, and speed-up, over the complete set of instances, per neighbourhood size. These
ablation runtimes are uncapped; the same-hardware comparison of \autoref{tab:comparison} caps runtimes at $120$s, so its All column at $n_g{=}24$ reads
slightly lower ($0.873$ vs.\ $0.935$).}
\label{tab:runtime_ng_stats}
\end{table}

\subsection{Hard Instances}
For a further drill down we consider the hard instances. Hard instances are defined as having a runtime of at least 10 seconds
 for the baseline algorithmic setting. 

The performance profile of the algorithm with the different configurations is shown in \ifijoc Figure~EC.6 of the online supplement\else\autoref{fig:performance-profile-hard}\fi.
It is evident that the combination of all optimisations are clearly superior to the other configurations on hard instances,
in all cases this is the best setting. It is also seen, that the vectorised dominance comparison seems to be more efficient
 on hard instances compared to the complete set of instances, although it is still the least efficient of the three optimisations.

\ifijoc\else
\begin{figure}[htbp]
    \centering
    \includegraphics[width=0.8\textwidth]{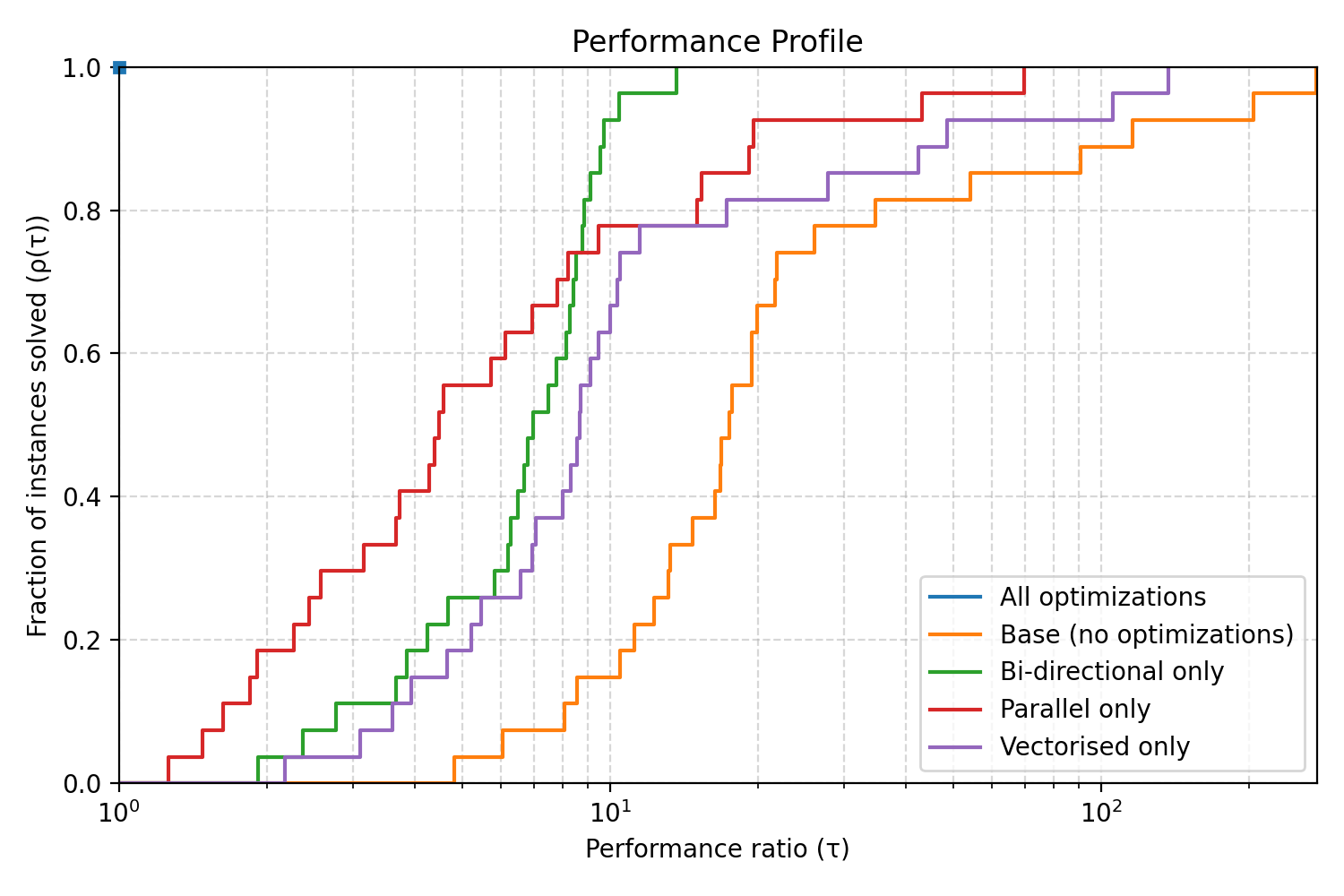}
    \caption{Performance profile of the algorithm with different configurations on hard instances}
    \label{fig:performance-profile-hard}
  \end{figure}
\fi

 \autoref{tab:runtime_ng_stats_hard} shows the runtime shifted (1s) geometric mean (s) unscaled, scaled, and speed-up for the hard instances per neighbourhood size.
The speed-ups on the hard instances are significantly higher than when considering the complete set of instances:
the combination of all optimisations is around 15--20x faster than the baseline (overall around 18x).

\begin{table}[htbp]
\centering
\begin{tabular}{lrrrrrr}
\toprule
NG & Base & Parallel & Bi-directional & Vectorised & All \\
\midrule
\multicolumn{6}{l}{Runtime (s)} \\
\midrule
8 & 31.779 & 6.206 & 10.624 & 14.243 & 2.130 \\
16 & 151.754 & 44.944 & 46.770 & 79.467 & 7.853 \\
24 & 328.995 & 92.481 & 94.210 & 157.949 & 16.444 \\
\midrule
\multicolumn{6}{l}{Normalized} \\
\midrule
8 & 14.920 & 2.914 & 4.988 & 6.687 & 1.000 \\
16 & 19.323 & 5.723 & 5.955 & 10.119 & 1.000 \\
24 & 20.007 & 5.624 & 5.729 & 9.605 & 1.000 \\
\midrule
\multicolumn{6}{l}{Speedup} \\
\midrule
8 & 1.000 & 5.120 & 2.991 & 2.231 & 14.920 \\
16 & 1.000 & 3.377 & 3.245 & 1.910 & 19.323 \\
24 & 1.000 & 3.557 & 3.492 & 2.083 & 20.007 \\
\bottomrule
\end{tabular}
\vspace{1em}
\caption{Runtime shifted (1s) geometric mean (s) unscaled, scaled, and speed-up, over the hard instances only (Base runtime $\geq 10$s), per neighbourhood size.}
\label{tab:runtime_ng_stats_hard}
\end{table}

\section{Conclusion}
\label{sec:conclusion}

We have presented a parallel pull labelling algorithm for the RCSPP built on an acyclic dependency bucket graph that makes label storage immutable,
together with a bi-directional extension whose midpoint emerges from the bucket processing order and a vectorised dominance criterion. On instances
derived from Solomon's VRPTW benchmarks, each optimisation gives around $1.4$--$2\times$ speed-up, the combination around $4\times$ on average and around
$18\times$ on the hard, long-time-window instances (reaching $274\times$ on the instance that benefits most, RC208). Against an open implementation of
the state-of-the-art bucket graph labelling algorithm on the same hardware, the pull algorithm is $1.9\times$--$2.4\times$ faster -- the expected, modest
margin of a specialised solver over a general library.

On the hardest large-neighbourhood instances the acyclic-dependency restriction limits parallelism and the algorithm shows a heavier worst-case tail;
closing it by permitting bounded cycles in the dependency graph, in the spirit of \citet{sadykov2020}, is the most consequential direction for future
work, though it risks the immutability and storage efficiency the scheme relies on. More generally, the exposed parallelism is governed by the shape of
the dependency graph: the forward, backward, and splice jobs of a vertex each form a chain, so at most a small constant number of jobs per vertex is in
progress at once, and the achievable speed-up is at most the total number of jobs divided by the longest chain of dependent jobs (the span). Broad,
shallow graphs therefore parallelise well and long, path-like ones poorly; smaller buckets refine granularity but lengthen the per-vertex chains rather
than raising this ceiling. Characterising the relationship precisely is open.

A further direction is the integration of non-robust cuts. Limited-memory rank-1 cuts \citep{pecin2017lmr1c}, which generalise the subset-row
inequalities of \citet{jepsen2008}, require reduced-cost-compensated dominance \citep{pecin2017}; our implementation already supports this through a
soft-dominance variant, but a study within a full branch-cut-and-price solver \citep{pessoa2020,costa2019}, including the vectorised handling of the
compensation, remains to be done. Finally, offloading the dominance check to the GPU may help for very large label sets, though memory-transfer overhead
may be a bottleneck.

\bibliographystyle{unsrt}
\bibliography{references}

@misc{bgspprc,
  author       = {Spoorendonk, Simon},
  title        = {{bucket-graph-spprc}: a header-only {C++23} bucket graph labeling library for the {SPPRC}},
  year         = {2026},
  publisher    = {Zenodo},
  version      = {v0.1.0},
  doi          = {10.5281/zenodo.20819208},
  note         = {Concept DOI (all versions); version v0.1.0 at \url{https://doi.org/10.5281/zenodo.20819209}.
                  Source at \url{https://github.com/spoorendonk/bucket-graph-spprc}, MIT license}
}

@misc{bgspprcpaper,
  author        = {Spoorendonk, Simon},
  title         = {\texttt{bucket-graph-spprc}: an extensible {C++} library for the shortest path problem with resource constraints},
  year          = {2026},
  eprint        = {2606.30847},
  archivePrefix = {arXiv},
  doi           = {10.48550/arXiv.2606.30847},
  howpublished  = {Preprint, arXiv:2606.30847},
  url           = {https://arxiv.org/abs/2606.30847}
}

@incollection{desrosiers1995,
  author    = {Desrosiers, Jacques and Dumas, Yvan and Solomon, Marius M. and Soumis, Fran{\c{c}}ois},
  title     = {Time Constrained Routing and Scheduling},
  booktitle = {Network Routing},
  series    = {Handbooks in Operations Research and Management Science},
  volume    = {8},
  publisher = {Elsevier},
  year      = {1995},
  pages     = {35--139},
  doi       = {10.1016/S0927-0507(05)80106-9}
}

@techreport{uchoa2024,
  title={{Optimizing with Column Generation}: Advanced Branch-Cut-and-Price Algorithms ({Part I})}, 
  author={Uchoa, Eduardo and Pessoa, Artur and Moreno, Lorenza}, 
  institution={Cadernos do LOGIS-UFF, Universidade Federal Fluminense, Engenharia de Produ{\c{c}}{\~a}o}, 
  number={L-2024-3}, 
  month={August}, 
  year={2024} 
}

@techreport{desrosiers2024,
  author = {Desrosiers, Jacques and Lübbecke, Marco and Desaulniers, Guy and Gauthier, Jean Bertrand},
  pages = {1-657},
  title = {Branch-and-Price},
  pubstate = {published},
  is_member_collaboration = {true},
  year = {2024},
  month = jun,
  month_numeric = {6},
  revision_year = {2024},
  revision_month = {Oct},
  edition = {Revised},
  number = {G-2024-36},
  type = {Les Cahiers du GERAD},
  institution = {Groupe d’études et de recherche en analyse des décisions},
  address = {GERAD, Montréal QC H3T 2A7, Canada},
  URL = {https://www.gerad.ca/en/papers/G-2024-36},
  eprint = {https://www.gerad.ca/papers/G-2024-36.pdf?locale=en}
}

@article{baldacci2011,
  title={New Route Relaxation and Pricing Strategies for the Vehicle Routing Problem},
  author={Baldacci, Roberto and Mingozzi, Aristide and Roberti, Roberto},
  journal={Operations Research},
  volume={59},
  number={5},
  pages={1269--1283},
  year={2011},
  doi={10.1287/opre.1110.0975},
  url={https://doi.org/10.1287/opre.1110.0975}
}

@article{pecin2017,
  title={Improved branch-cut-and-price for capacitated vehicle routing},
  author={Pecin, Diego and Pessoa, Artur and Poggi, Marcus and Uchoa, Eduardo},
  journal={Mathematical Programming Computation},
  volume={9},
  number={1},
  pages={61--100},
  year={2017},
  doi={10.1007/s12532-016-0108-8},
  publisher={Springer}
}

@article{jepsen2008,
  title={Subset-Row Inequalities Applied to the Vehicle-Routing Problem with Time Windows},
  author={Jepsen, Mads and Petersen, Bjørn and Spoorendonk, Simon and Pisinger, David},
  journal={Operations Research},
  volume={56},
  number={2},
  pages={497--511},
  year={2008},
  doi={10.1287/opre.1070.0449},
  publisher={INFORMS}
}

@article{pathwyse,
  title={{PathWyse}: a flexible, open-source library for the resource constrained shortest path problem},
  author={Salani, Matteo and Basso, Saverio and Giuffrida, Vincenzo},
  journal={Optimization Methods and Software},
  volume={39},
  number={2},
  pages={298--320},
  year={2024},
  doi={10.1080/10556788.2023.2296978},
  publisher={Taylor \& Francis}
}

@article{lu2021,
  title={Accelerating Exact Constrained Shortest Paths on {GPUs}},
  author={Lu, Shengliang and He, Bingsheng and Li, Yuchen and Fu, Hao},
  journal={Proceedings of the VLDB Endowment},
  volume={14},
  number={4},
  pages={547--559},
  year={2021},
  doi={10.14778/3436905.3436914},
  publisher={VLDB Endowment}
}

@article{sadykov2020,
  title={A Bucket Graph--Based Labeling Algorithm with Application to Vehicle Routing},
  author={Sadykov, Ruslan and Uchoa, Eduardo and Pessoa, Artur},
  journal={Transportation Science},
  volume={55},
  number={1},
  pages={4--28},
  year={2021},
  doi={10.1287/trsc.2020.0985},
  publisher={INFORMS}
}

@article{pecin2017lmr1c,
  title={Limited memory rank-1 cuts for vehicle routing problems},
  author={Pecin, Diego and Pessoa, Artur and Poggi, Marcus and Uchoa, Eduardo and Santos, Haroldo},
  journal={Operations Research Letters},
  volume={45},
  number={3},
  pages={206--209},
  year={2017},
  doi={10.1016/j.orl.2017.02.006},
  publisher={Elsevier}
}

@article{pessoa2020,
  title={A Generic Exact Solver for Vehicle Routing and Related Problems},
  author={Pessoa, Artur and Sadykov, Ruslan and Uchoa, Eduardo and Vanderbeck, Fran{\c{c}}ois},
  journal={Mathematical Programming},
  volume={183},
  number={1-2},
  pages={483--523},
  year={2020},
  doi={10.1007/s10107-020-01523-z},
  publisher={Springer}
}

@article{costa2019,
  title={Exact Branch-Price-and-Cut Algorithms for Vehicle Routing},
  author={Costa, Luciano and Contardo, Claudio and Desaulniers, Guy},
  journal={Transportation Science},
  volume={53},
  number={4},
  pages={946--985},
  year={2019},
  doi={10.1287/trsc.2018.0878},
  publisher={INFORMS}
}

@article{pecin2017b,
  title={New Enhancements for the Exact Solution of the Vehicle Routing Problem with Time Windows},
  author={Pecin, Diego and Contardo, Claudio and Desaulniers, Guy and Uchoa, Eduardo},
  journal={INFORMS Journal on Computing},
  volume={29},
  number={3},
  pages={489--502},
  year={2017},
  doi={10.1287/ijoc.2016.0744},
  publisher={INFORMS}
}

@incollection{irnich2005,
  author    = {Stefan Irnich and Guy Desaulniers},
  title     = {Shortest Path Problems with Resource Constraints},
  booktitle = {Column Generation},
  editor    = {Guy Desaulniers and Jacques Desrosiers and Marius M. Solomon},
  year      = {2005},
  publisher = {Springer US},
  address   = {Boston, MA},
  pages     = {33--65},
  isbn      = {978-0-387-25486-9},
  doi       = {10.1007/0-387-25486-2_2},
  url       = {https://doi.org/10.1007/0-387-25486-2_2}
}

@book{desaulniers2005,
  title        = {Column Generation},
  editor       = {Guy Desaulniers and Jacques Desrosiers and Marius M. Solomon},
  year         = {2005},
  publisher    = {Springer},
  address      = {Boston, MA},
  doi          = {10.1007/b135457},
  isbn         = {978-0-387-25486-9},
  series       = {Business and Economics, Business and Management (R0)},
  edition      = {1},
  pages        = {XV + 358},
  url          = {https://doi.org/10.1007/b135457}
}

@article{feillet2004,
author = {Feillet, Dominique and Dejax, Pierre and Gendreau, Michel and Gueguen, Cyrille},
title = {An exact algorithm for the elementary shortest path problem with resource constraints: Application to some vehicle routing problems},
journal = {Networks},
volume = {44},
number = {3},
pages = {216-229},
doi = {10.1002/net.20033},
url = {https://onlinelibrary.wiley.com/doi/abs/10.1002/net.20033},
year = {2004}
}

@article{salani2006,
title = {Symmetry helps: Bounded bi-directional dynamic programming for the elementary shortest path problem with resource constraints},
journal = {Discrete Optimization},
volume = {3},
number = {3},
pages = {255-273},
year = {2006},
note = {Graphs and Combinatorial Optimization},
issn = {1572-5286},
doi = {10.1016/j.disopt.2006.05.007},
url = {https://www.sciencedirect.com/science/article/pii/S1572528606000417},
author = {Giovanni Righini and Matteo Salani}
}

@article{tilk2017,
title = {Asymmetry matters: Dynamic half-way points in bidirectional labeling for solving shortest path problems with resource constraints faster},
journal = {European Journal of Operational Research},
volume = {261},
number = {2},
pages = {530-539},
year = {2017},
issn = {0377-2217},
doi = {10.1016/j.ejor.2017.03.017},
url = {https://www.sciencedirect.com/science/article/pii/S0377221717302035},
author = {Christian Tilk and Ann-Kathrin Rothenbächer and Timo Gschwind and Stefan Irnich}
}

@article{salani2024,
title = {Enhanced bi-directional dynamic programming algorithm for the resource constrained shortest path problem},
journal = {Transportation Research Procedia},
volume = {78},
pages = {361-368},
year = {2024},
note = {25th Euro Working Group on Transportation Meeting},
issn = {2352-1465},
doi = {10.1016/j.trpro.2024.02.046},
url = {https://www.sciencedirect.com/science/article/pii/S2352146524000991},
author = {Matteo Salani and Saverio Basso and Giovanni Righini}
}

@article{bulhoes2018,
title = {A branch-and-price algorithm for the Minimum Latency Problem},
journal = {Computers \& Operations Research},
volume = {93},
pages = {66-78},
year = {2018},
issn = {0305-0548},
doi = {10.1016/j.cor.2018.01.016},
url = {https://www.sciencedirect.com/science/article/pii/S0305054818300248},
author = {Teobaldo Bulhões and Ruslan Sadykov and Eduardo Uchoa}
}

@article{irnich2008,
  author    = {Stefan Irnich},
  title     = {Resource extension functions: properties, inversion, and generalization to segments},
  journal   = {OR Spectrum},
  year      = {2008},
  volume    = {30},
  number    = {1},
  pages     = {113--148},
  doi       = {10.1007/s00291-007-0083-6},
  url       = {https://doi.org/10.1007/s00291-007-0083-6},
  issn      = {1436-6304}
}

@article{solomon1987,
  title={Algorithms for the Vehicle Routing and Scheduling Problems with Time Window Constraints},
  author={Solomon, Marius M.},
  journal={Operations Research},
  volume={35},
  number={2},
  pages={254--265},
  year={1987},
  doi={10.1287/opre.35.2.254},
  publisher={INFORMS}
}

@misc{spoorendonk2025,
  author       = {Spoorendonk, Simon},
  title        = {Resource Constrained Shortest Path Problem Instances with Time Windows, Capacity, and Neighbourhoods},
  year         = {2025},
  howpublished = {\url{https://github.com/spoorendonk/rcspp_dataset}},
  note         = {Dataset}
}

\appendix
\section{Appendix}

\subsection{Runtime Details}
\label{sec:runtime-details}

\autoref{tab:runtime_stats_8}, \autoref{tab:runtime_stats_16}, and \autoref{tab:runtime_stats_24} show the detailed
runtime details for the instances.

\begin{longtable}{lrrrrrr}
\toprule
Instance & Base & Vectorised & Parallel & Bi-directional & All \\
\midrule
\endfirsthead
\toprule
Instance & Base & Vectorised & Parallel & Bi-directional & All \\
\midrule
\endhead
\bottomrule
\addlinespace
\multicolumn{7}{r}{{Continued on next page}} \\
\endfoot
\bottomrule
\addlinespace
\caption{Runtime (s) by instance, maximum neighbourhood size N8.}
\label{tab:runtime_stats_8}
\endlastfoot
C101 & 0.008 & 0.007 & 0.012 & 0.005 & 0.010 \\
C102 & 0.036 & 0.030 & 0.028 & 0.015 & 0.017 \\
C103 & 0.072 & 0.059 & 0.037 & 0.039 & 0.023 \\
C104 & 0.377 & 0.249 & 0.128 & 0.171 & 0.057 \\
C105 & 0.008 & 0.008 & 0.011 & 0.004 & 0.009 \\
C106 & 0.012 & 0.011 & 0.013 & 0.007 & 0.010 \\
C107 & 0.010 & 0.010 & 0.012 & 0.006 & 0.009 \\
C108 & 0.027 & 0.023 & 0.021 & 0.013 & 0.013 \\
C109 & 0.056 & 0.046 & 0.030 & 0.028 & 0.017 \\
\midrule
C201 & 0.008 & 0.008 & 0.015 & 0.005 & 0.010 \\
C202 & 0.111 & 0.081 & 0.079 & 0.048 & 0.030 \\
C203 & 1.713 & 0.881 & 0.932 & 0.666 & 0.163 \\
C204 & 36.388 & 16.169 & 10.882 & 10.844 & 2.961 \\
C205 & 0.072 & 0.055 & 0.049 & 0.028 & 0.020 \\
C206 & 0.147 & 0.097 & 0.091 & 0.066 & 0.037 \\
C207 & 0.491 & 0.264 & 0.304 & 0.309 & 0.211 \\
C208 & 0.233 & 0.151 & 0.139 & 0.093 & 0.044 \\
\midrule
R101 & 0.003 & 0.003 & 0.007 & 0.002 & 0.006 \\
R102 & 0.015 & 0.014 & 0.014 & 0.009 & 0.011 \\
R103 & 0.045 & 0.040 & 0.024 & 0.024 & 0.019 \\
R104 & 0.532 & 0.314 & 0.088 & 0.149 & 0.038 \\
R105 & 0.010 & 0.009 & 0.011 & 0.005 & 0.008 \\
R106 & 0.046 & 0.037 & 0.024 & 0.021 & 0.017 \\
R107 & 0.199 & 0.139 & 0.057 & 0.091 & 0.037 \\
R108 & 1.009 & 0.558 & 0.147 & 0.276 & 0.063 \\
R109 & 0.052 & 0.035 & 0.024 & 0.014 & 0.013 \\
R110 & 0.232 & 0.136 & 0.054 & 0.038 & 0.020 \\
R111 & 0.204 & 0.135 & 0.054 & 0.079 & 0.037 \\
R112 & 0.851 & 0.460 & 0.131 & 0.176 & 0.047 \\
\midrule
R201 & 0.026 & 0.022 & 0.025 & 0.015 & 0.016 \\
R202 & 0.201 & 0.141 & 0.096 & 0.177 & 0.084 \\
R203 & 0.921 & 0.548 & 0.337 & 0.720 & 0.293 \\
R204 & 11.922 & 5.486 & 2.270 & 8.104 & 1.394 \\
R205 & 0.326 & 0.203 & 0.143 & 0.222 & 0.070 \\
R206 & 0.763 & 0.457 & 0.262 & 0.680 & 0.184 \\
R207 & 5.078 & 2.465 & 1.157 & 4.037 & 0.754 \\
R208 & 22.232 & 9.891 & 3.146 & 13.815 & 2.124 \\
R209 & 0.752 & 0.422 & 0.204 & 0.595 & 0.130 \\
R210 & 0.883 & 0.506 & 0.267 & 0.692 & 0.178 \\
R211 & 11.782 & 5.198 & 1.428 & 2.801 & 0.452 \\
\midrule
RC101 & 0.008 & 0.008 & 0.009 & 0.004 & 0.008 \\
RC102 & 0.043 & 0.033 & 0.021 & 0.015 & 0.014 \\
RC103 & 0.223 & 0.137 & 0.050 & 0.051 & 0.024 \\
RC104 & 3.720 & 1.742 & 0.733 & 0.440 & 0.075 \\
RC105 & 0.024 & 0.020 & 0.018 & 0.010 & 0.012 \\
RC106 & 0.054 & 0.040 & 0.021 & 0.014 & 0.013 \\
RC107 & 0.787 & 0.378 & 0.174 & 0.054 & 0.022 \\
RC108 & 2.042 & 0.959 & 0.481 & 0.167 & 0.059 \\
\midrule
RC201 & 0.028 & 0.026 & 0.024 & 0.019 & 0.018 \\
RC202 & 0.238 & 0.167 & 0.101 & 0.179 & 0.062 \\
RC203 & 5.138 & 2.407 & 1.960 & 2.840 & 0.744 \\
RC204 & 88.486 & 39.552 & 20.290 & 46.638 & 10.969 \\
RC205 & 0.131 & 0.096 & 0.064 & 0.078 & 0.038 \\
RC206 & 0.276 & 0.182 & 0.110 & 0.176 & 0.058 \\
RC207 & 1.755 & 0.898 & 0.493 & 0.841 & 0.162 \\
RC208 & 95.616 & 40.142 & 15.820 & 7.528 & 0.826 \\
\end{longtable}

\newpage
\begin{longtable}{lrrrrrr}
\toprule
Instance & Base & Vectorised & Parallel & Bi-directional & All \\
\midrule
\endfirsthead
\toprule
Instance & Base & Vectorised & Parallel & Bi-directional & All \\
\midrule
\endhead
\bottomrule
\addlinespace
\multicolumn{7}{r}{{Continued on next page}} \\
\endfoot
\bottomrule
\addlinespace
\caption{Runtime (s) by instance, maximum neighbourhood size N16.}
\label{tab:runtime_stats_16}
\endlastfoot
C101 & 0.007 & 0.011 & 0.012 & 0.005 & 0.010 \\
C102 & 0.030 & 0.027 & 0.023 & 0.015 & 0.016 \\
C103 & 0.077 & 0.064 & 0.039 & 0.042 & 0.026 \\
C104 & 2.676 & 1.582 & 1.561 & 0.730 & 0.319 \\
C105 & 0.008 & 0.008 & 0.012 & 0.004 & 0.009 \\
C106 & 0.011 & 0.011 & 0.014 & 0.007 & 0.010 \\
C107 & 0.011 & 0.010 & 0.013 & 0.006 & 0.009 \\
C108 & 0.026 & 0.023 & 0.021 & 0.013 & 0.012 \\
C109 & 0.057 & 0.048 & 0.031 & 0.030 & 0.017 \\
\midrule
C201 & 0.008 & 0.008 & 0.015 & 0.005 & 0.010 \\
C202 & 0.169 & 0.122 & 0.141 & 0.062 & 0.040 \\
C203 & 48.067 & 25.596 & 38.067 & 5.872 & 2.478 \\
C204 & 274.322 & 139.452 & 93.462 & 76.135 & 16.294 \\
C205 & 0.065 & 0.053 & 0.051 & 0.027 & 0.020 \\
C206 & 0.151 & 0.115 & 0.108 & 0.051 & 0.033 \\
C207 & 0.504 & 0.322 & 0.504 & 0.544 & 0.277 \\
C208 & 0.151 & 0.114 & 0.101 & 0.072 & 0.040 \\
\midrule
R101 & 0.003 & 0.003 & 0.007 & 0.002 & 0.006 \\
R102 & 0.014 & 0.014 & 0.015 & 0.009 & 0.012 \\
R103 & 0.056 & 0.050 & 0.028 & 0.024 & 0.019 \\
R104 & 1.209 & 0.769 & 0.248 & 0.313 & 0.063 \\
R105 & 0.010 & 0.010 & 0.011 & 0.005 & 0.008 \\
R106 & 0.052 & 0.045 & 0.028 & 0.022 & 0.017 \\
R107 & 0.182 & 0.145 & 0.063 & 0.118 & 0.042 \\
R108 & 1.037 & 0.678 & 0.185 & 0.333 & 0.073 \\
R109 & 0.042 & 0.035 & 0.022 & 0.015 & 0.013 \\
R110 & 0.193 & 0.139 & 0.055 & 0.040 & 0.020 \\
R111 & 0.212 & 0.159 & 0.068 & 0.097 & 0.041 \\
R112 & 0.881 & 0.559 & 0.152 & 0.228 & 0.062 \\
\midrule
R201 & 0.025 & 0.024 & 0.024 & 0.016 & 0.017 \\
R202 & 0.271 & 0.217 & 0.135 & 0.254 & 0.101 \\
R203 & 2.856 & 1.748 & 0.937 & 2.402 & 0.816 \\
R204 & 76.611 & 40.083 & 14.119 & 43.329 & 5.782 \\
R205 & 0.515 & 0.334 & 0.219 & 0.364 & 0.105 \\
R206 & 2.844 & 1.705 & 0.892 & 3.006 & 0.589 \\
R207 & 38.786 & 20.676 & 10.663 & 20.226 & 2.372 \\
R208 & 401.236 & 200.622 & 78.720 & 257.228 & 30.547 \\
R209 & 2.649 & 1.597 & 1.097 & 1.689 & 0.274 \\
R210 & 5.092 & 2.982 & 1.647 & 3.636 & 0.515 \\
R211 & 100.435 & 51.327 & 17.512 & 16.440 & 1.852 \\
\midrule
RC101 & 0.008 & 0.007 & 0.009 & 0.004 & 0.008 \\
RC102 & 0.041 & 0.034 & 0.023 & 0.015 & 0.014 \\
RC103 & 0.265 & 0.190 & 0.066 & 0.059 & 0.026 \\
RC104 & 9.194 & 5.000 & 3.664 & 0.859 & 0.262 \\
RC105 & 0.023 & 0.021 & 0.019 & 0.010 & 0.014 \\
RC106 & 0.050 & 0.041 & 0.021 & 0.014 & 0.013 \\
RC107 & 0.530 & 0.333 & 0.151 & 0.056 & 0.023 \\
RC108 & 2.266 & 1.272 & 0.759 & 0.185 & 0.059 \\
\midrule
RC201 & 0.027 & 0.027 & 0.027 & 0.021 & 0.019 \\
RC202 & 0.316 & 0.220 & 0.134 & 0.147 & 0.061 \\
RC203 & 15.521 & 8.427 & 6.424 & 6.181 & 0.924 \\
RC204 & 538.427 & 275.843 & 201.953 & 620.551 & 88.971 \\
RC205 & 0.280 & 0.200 & 0.156 & 0.117 & 0.047 \\
RC206 & 0.477 & 0.319 & 0.203 & 0.229 & 0.063 \\
RC207 & 5.479 & 3.134 & 1.670 & 3.047 & 0.515 \\
RC208 & 2984.494 & 1543.195 & 629.381 & 128.203 & 14.600 \\
\end{longtable}

\begin{longtable}{lrrrrrr}
\toprule
Instance & Base & Vectorised & Parallel & Bi-directional & All \\
\midrule
\endfirsthead
\toprule
Instance & Base & Vectorised & Parallel & Bi-directional & All \\
\midrule
\endhead
\bottomrule
\addlinespace
\multicolumn{7}{r}{{Continued on next page}} \\
\endfoot
\bottomrule
\addlinespace
\caption{Runtime (s) by instance, maximum neighbourhood size N24.}
\label{tab:runtime_stats_24}
\endlastfoot
C101 & 0.007 & 0.007 & 0.011 & 0.004 & 0.010 \\
C102 & 0.020 & 0.020 & 0.019 & 0.013 & 0.017 \\
C103 & 0.162 & 0.122 & 0.071 & 0.056 & 0.030 \\
C104 & 0.971 & 0.560 & 0.624 & 0.321 & 0.121 \\
C105 & 0.008 & 0.008 & 0.015 & 0.004 & 0.009 \\
C106 & 0.012 & 0.012 & 0.014 & 0.006 & 0.010 \\
C107 & 0.010 & 0.011 & 0.012 & 0.006 & 0.010 \\
C108 & 0.031 & 0.023 & 0.021 & 0.013 & 0.013 \\
C109 & 0.051 & 0.042 & 0.029 & 0.030 & 0.017 \\
\midrule
C201 & 0.007 & 0.008 & 0.015 & 0.005 & 0.010 \\
C202 & 0.783 & 0.437 & 0.895 & 0.186 & 0.084 \\
C203 & 11.393 & 5.413 & 8.575 & 1.096 & 0.571 \\
C204 & 3498.628 & 1678.476 & 1947.330 & 915.003 & 237.477 \\
C205 & 0.047 & 0.038 & 0.038 & 0.020 & 0.017 \\
C206 & 0.113 & 0.093 & 0.083 & 0.048 & 0.030 \\
C207 & 0.176 & 0.130 & 0.125 & 0.156 & 0.095 \\
C208 & 0.268 & 0.175 & 0.159 & 0.083 & 0.041 \\
\midrule
R101 & 0.003 & 0.003 & 0.006 & 0.002 & 0.006 \\
R102 & 0.014 & 0.014 & 0.014 & 0.009 & 0.012 \\
R103 & 0.061 & 0.053 & 0.028 & 0.026 & 0.018 \\
R104 & 1.155 & 0.672 & 0.191 & 0.289 & 0.056 \\
R105 & 0.033 & 0.010 & 0.011 & 0.005 & 0.008 \\
R106 & 0.046 & 0.041 & 0.026 & 0.021 & 0.018 \\
R107 & 0.184 & 0.151 & 0.056 & 0.096 & 0.038 \\
R108 & 1.805 & 0.990 & 0.318 & 0.442 & 0.081 \\
R109 & 0.043 & 0.036 & 0.025 & 0.014 & 0.013 \\
R110 & 0.218 & 0.152 & 0.058 & 0.041 & 0.020 \\
R111 & 0.235 & 0.164 & 0.065 & 0.094 & 0.042 \\
R112 & 0.880 & 0.513 & 0.149 & 0.212 & 0.063 \\
\midrule
R201 & 0.024 & 0.023 & 0.023 & 0.017 & 0.017 \\
R202 & 0.566 & 0.382 & 0.237 & 0.481 & 0.139 \\
R203 & 38.077 & 17.190 & 9.952 & 21.806 & 7.891 \\
R204 & 1369.055 & 637.270 & 233.351 & 766.177 & 121.985 \\
R205 & 0.479 & 0.310 & 0.192 & 0.364 & 0.092 \\
R206 & 37.330 & 17.063 & 7.831 & 23.331 & 1.706 \\
R207 & 233.631 & 104.480 & 44.873 & 125.355 & 12.027 \\
R208 & 12959.341 & 6462.977 & 2914.763 & 2908.629 & 373.768 \\
R209 & 4.463 & 2.247 & 1.043 & 2.151 & 0.296 \\
R210 & 24.337 & 11.138 & 5.971 & 11.518 & 1.392 \\
R211 & 1016.520 & 474.254 & 218.907 & 91.289 & 11.193 \\
\midrule
RC101 & 0.008 & 0.007 & 0.011 & 0.004 & 0.008 \\
RC102 & 0.041 & 0.033 & 0.021 & 0.015 & 0.016 \\
RC103 & 0.212 & 0.143 & 0.048 & 0.053 & 0.025 \\
RC104 & 9.282 & 4.366 & 2.900 & 0.847 & 0.230 \\
RC105 & 0.026 & 0.022 & 0.016 & 0.010 & 0.014 \\
RC106 & 0.053 & 0.040 & 0.021 & 0.013 & 0.013 \\
RC107 & 0.539 & 0.295 & 0.137 & 0.050 & 0.022 \\
RC108 & 3.499 & 1.706 & 1.069 & 0.197 & 0.060 \\
\midrule
RC201 & 0.033 & 0.029 & 0.027 & 0.021 & 0.018 \\
RC202 & 0.216 & 0.163 & 0.092 & 0.149 & 0.059 \\
RC203 & 9.105 & 4.308 & 2.142 & 4.976 & 1.158 \\
RC204 & 1667.897 & 807.095 & 471.822 & 737.057 & 77.071 \\
RC205 & 0.161 & 0.121 & 0.075 & 0.080 & 0.038 \\
RC206 & 0.536 & 0.338 & 0.171 & 0.263 & 0.071 \\
RC207 & 13.884 & 6.532 & 3.440 & 5.328 & 0.784 \\
RC208 & 9620.517 & 4808.602 & 2445.809 & 341.233 & 35.112 \\
\end{longtable}

\end{document}